\documentclass[reqno,11pt]{amsart}

\usepackage{amsfonts,amsmath,amsthm,amssymb}
\usepackage{shortvrb,graphicx,psfrag}
\usepackage{wrapfig}
\usepackage[colorlinks,pdfpagelabels,pdfstartview = FitH,bookmarksopen
= true,bookmarksnumbered = true,linkcolor = blue,plainpages =
false,hypertexnames = false,citecolor = red,pagebackref=false]{hyperref}
\usepackage{relsize}

\usepackage{xcolor}
\usepackage{appendix}
\usepackage{comment}
\usepackage{ulem}
\usepackage[makeroom]{cancel}
\usepackage[margin=1.4in]{geometry}
\usepackage{stmaryrd}
\usepackage{float}
\usepackage{mathrsfs}
\SetSymbolFont{stmry}{bold}{U}{stmry}{m}{n}

\allowdisplaybreaks

\let\TeXchi\chi
\newbox\chibox
\setbox0 \hbox{\mathsurround0pt $\TeXchi$}
\setbox\chibox \hbox{\raise\dp0 \box 0 }
\def\chi{\copy\chibox}

\def\Xiint#1{\mathchoice
    {\XXiint\displaystyle\textstyle{#1}}%
    {\XXiint\textstyle\scriptstyle{#1}}%
    {\XXiint\scriptstyle\scriptscriptstyle{#1}}%
    {\XXiint\scriptscriptstyle\scriptscriptstyle{#1}}%
    \!\iint}
\def\XXiint#1#2#3{\setbox0=\hbox{$#1{#2#3}{\iint}$}
    \vcenter{\hbox{$#2#3$}}\kern-0.5\wd0}
\def\biint{\Xiint{-\!-}}

\renewcommand{\d}{\mathrm{d}}
\newcommand{\dx}{\mathrm{d}x}
\newcommand{\dy}{\mathrm{d}y}

\newcommand{\dt}{\mathrm{d}t}

\renewcommand{\rho}{\varrho}

\newcommand{\noi}{\noindent}
\newcommand{\dsty}{\displaystyle}

\newcommand{\pl}{\partial}

\newcommand{\al}{\alpha}

\newcommand{\be}{\beta}

\newcommand{\gm}{\gamma}
\newcommand{\dl}{\delta}
\newcommand{\Dl}{\Delta}
\newcommand{\lm}{\lambda}

\newcommand{\varep}{\varepsilon}

\newcommand{\vp}{\varphi}
\newcommand{\sig}{\sigma}

\newcommand{\om}{\omega}
\newcommand{\Om}{\Omega}

\newcommand{\z}{\zeta}

\newcommand{\loc}{\operatorname{loc}}

\newcommand{\dvg}{\operatorname{div}}
\newcommand{\essosc}{\operatornamewithlimits{ess\,osc}}
\newcommand{\essup}{\operatornamewithlimits{ess\,sup}}
\newcommand{\essinf}{\operatornamewithlimits{ess\,inf}}

\newcommand{\rr}{\mathbb{R}}
\newcommand{\rn}{\rr^N}
\newcommand{\nn}{\mathbb{N}}

\def\Xint#1{\mathchoice
    {\XXint\displaystyle\textstyle{#1}}
    {\XXint\textstyle\scriptstyle{#1}}
    {\XXint\scriptstyle\scriptscriptstyle{#1}}
    {\XXint\scriptscriptstyle\scriptscriptstyle{#1}}
    \!\int}
\def\XXint#1#2#3{\setbox0=\hbox{$#1{#2#3}{\int}$}
    \vcenter{\hbox{$#2#3$}}\kern-0.5\wd0}
\def\bint{\Xint-}
\def\dashint{\Xint{\raise4pt\hbox to7pt{\hrulefill}}}
\def\dashiint{\bint\kern-0.15cm\bint}

\newtheorem{proposition}{Proposition}[section]
\newtheorem{theorem}{Theorem}[section]

\newtheorem{lemma}{Lemma}[section]

\newtheorem{remark}{Remark}[section]
\newtheorem{definition}{Definition}[section]

\numberwithin{equation}{section}
\numberwithin{theorem}{section}
\numberwithin{proposition}{section}
\numberwithin{lemma}{section}
\numberwithin{remark}{section}
\title[Nonlocal Harnack estimates]{Harnack estimates for \\ nonlocal drift-diffusion equations}

\author[N. Liao]{Naian Liao}
\address{Fachbereich Mathematik, Paris-Lodron-Universit\"at Salzburg,
Hellbrunner Str. 34, 5020 Salzburg, Austria}{}
\email{naian.liao@plus.ac.at}

\begin{document}



\begin{abstract}
A set of pointwise estimates are established for local solutions to nonlocal diffusion equations with a drift term. 
In particular, our Harnack estimates are the first ones for such equations, and
our H\"older regularity refines certain known result in several aspects. 
The approach is measure theoretical in the spirit of DeGiorgi classes. 
It yields novel nonlocal weak Harnack estimates in the elliptic case as well.
\end{abstract}  

\date{\today}

\maketitle
\section{Introduction}
Let $E\subset\rn$ be a bounded open set and $E_T:=E\times(0,T]$ for some $T>0$. 
Consider the nonlocal diffusion equation with a drift term:
\begin{equation}\label{Eq:1:1}
\left\{
\begin{array}{cc}
\displaystyle\pl_t u(x,t) + {\rm p.v.}\int_{\rn} a(x,y,t)  \frac{u(x,t) - u(y,t)}{|x-y|^{N+2s}}\,\dy+b\cdot \nabla u(x,t)=0\quad\text{in}\> E_T,\\
\\
\dvg   b(\cdot, t)=0\quad\text{in}\>E,\quad\text{for any}\> t\in(0,T],
\end{array}\right.
\end{equation}
where $s\in(0,1)$. The coefficient $a(x,y,t)$ is a measurable function in $\rr^{2N}\times(0,T]$ that is  symmetric about $(x,y)$ and satisfies $C_o\le a\le C_1$ for some positive $C_o,\,C_1$, whereas the vector valued function $b\in L^{\infty}(0,T; L_{\loc}^{N/s}(E))$ satisfies for some $C_2>0$ that  
\begin{equation}\label{Eq:b:1}
\essup_{t\in[0,T]} \Big(\bint_{K_{R}(x_o)} |b(x,t)|^{\frac{N}s}\,\dx\Big)^{\frac{s}N}\le C_2 R^{1-2s}
\end{equation}
for any ball $K_R(x_o)\subset E$. 
Our notion of weak solution is a local and variational concept that is independent of any initial or boundary data; see \S~\ref{S:notion}. 

Consideration of equation~\eqref{Eq:1:1} is motivated by the so-called surface quasi-geostrophic model--a toy model for the regularity theory of the Navier-Stokes equations, cf.~\cite{CMT-94, Const-06}. For recent breakthrough on the regularity of solutions, see~\cite{Caff-Vass-10, KN-09, KNV-07}.
In this paper, we aim to establish first Harnack estimates for solutions to the nonlocal drift-diffusion equation \eqref{Eq:1:1}. The new argument is flexible and refines previous H\"older regularity results as well. In addition, even for the elliptic case it readily yields novel nonlocal weak Harnack estimates; see~\S~\ref{S:elliptic}. As is well known that Harnack estimates are stronger than H\"older regularity and depict more precisely the pointwise behavior of solutions. The Harnack estimates derived in this paper can be used to study initial traces, to characterize global behavior of solutions, and to obtain sub-potential lower bounds of fundamental solutions. 

Titled ``Nonlocal weak Harnack estimates," the first version of this article presents the drift-free case and can be found on arXiv.

\subsection{Weak Harnack estimates}
Now we discuss weak Harnack estimates for non-negative super-solutions. 
Note that we will state the results in a unified fashion for all $s\in(0,1)$ while keeping in mind that
when $s\in(0,\frac12)$, condition~\eqref{Eq:b:1} forces $b=0$; when $s=\frac12$, it is equivalent to $b\in L_{\loc}^{\infty}(E_{T})$; and when $s\in(\frac12,1)$, it requires $b(\cdot,t)$ to be in a local Morrey space that includes $L_{\loc}^{\frac{N}{2s-1}}(E)$ as a proper subset.
The first result is as follows.

\begin{theorem}\label{Prop:WHI:1}
Let $u$ be a local, weak super-solution to \eqref{Eq:1:1}   in $E_T$, such that $u\ge0$. Assume that $b$ satisfies \eqref{Eq:b:1}.
There exist  constants $\varep,\,\eta\in(0,1)$ depending on the data $\{s, N, C_o, C_1\}$ and also on $C_2$, 
such that 
\begin{equation*}
	u(x,t) \ge \eta\Big(\bint_{K_\rho(x_o)} u^{\varep}(x,t_o)\,\dx\Big)^{\frac1{\varep}}
\end{equation*}
for almost all
$$ 
(x,t)\in K_{2\varrho}(x_o) \times\big( t_o+\tfrac12 (8\varrho)^{2s}, t_o+2(8\varrho)^{2s}\big],
$$
provided  that
\[
K_{4\rho}(x_o)\times\big(t_o, t_o+4(8\varrho)^{2s}\big]\subset  E_T.
\]
\end{theorem}
\noi The above theorem is new even for $b=0$ in the sense that we do not need any information of solutions below the time level $t_o$. Such kind of result, even for local diffusion equations, is not part of the repertoire in Moser's approach. Very recently, we can improve the integral power in the drift-free case (see~\cite[Theorem~1.4]{LW-24}), but the method requires to use the equation more than once. This type of formulation of weak Harnack estimates is useful in studying initial traces of global solutions, cf.~\cite[Chap. XI, \S~4]{DB}.

The second weak Harnack estimate improves the integral power of the previous one, yet an additional time integration is imposed, and moreover, it requires information of solutions below the time integration.
\begin{theorem}\label{Prop:WHI:2-}
Let $u$ be a  local, weak super-solution to \eqref{Eq:1:1}   in $E_T$, such that $u\ge0$. Assume that $b$ satisfies \eqref{Eq:b:1}.
There exists a constant   
 $\eta\in(0,1)$ depending on the data $\{s, N, C_o, C_1\}$  and also on $C_2$, such that 
\begin{equation*}
	u(x,t)\ge\eta \biint_{(x_o,t_o)+Q_\rho} u\,\dx\dt
\end{equation*}
for almost all 
$$ 
(x,t)\in K_{\varrho}(x_o) \times\big( t_o+\tfrac34 (4\varrho)^{2s},
	t_o+(4\varrho)^{2s}\big],
$$
provided that  
\[
K_{2\rho}(x_o)\times\big(t_o-(2\varrho)^{2s}, t_o+6(4\varrho)^{2s}\big]\subset E_T.
\]
\end{theorem}
\noi When $s\in[\frac12,1)$, the above theorem generalizes an important result of Felsinger \& Kassmann to a drift-diffusion case; see \cite[Theorem~1.1]{Kass-13} which assumes some more general structures. Although their result also included a bounded forcing term $f$, it can also be readily treated by our approach. In fact, simple calculations show that it suffices to assume a certain integrability of $f$. This remark applies to all results in this paper.

Both of the previous two weak Harnack estimates examine the gain of positivity from local integrals of the super-solution,
whereas the third weak Harnack estimate reveals an additional positivity contribution from the super-solution's long-range behavior; see~\eqref{Eq:tail} for the definition of tail.
As such, it presents a strong nonlocal character.
\begin{theorem}\label{Prop:WHI:2}
Let $u$ be a  local, weak super-solution to \eqref{Eq:1:1}   in $E_T$, such that $u\ge0$. Assume that $b$ satisfies \eqref{Eq:b:1}.
There exists a constant   
 $\eta\in(0,1)$ depending on the data $\{s, N, C_o, C_1\}$ and also on $C_2$, such that 
\begin{equation*}
	u(x,t)\ge\eta {\rm Tail}[u; (x_o,t_o)+Q_{\rho}]
\end{equation*}
for almost all
$$
(x,t)\in K_{\varrho}(x_o) \times\big( t_o+\tfrac34 (4\varrho)^{2s},
	t_o+(4\varrho)^{2s}\big],
$$
provided that  
\[
K_{2\rho}(x_o)\times\big(t_o-(2\varrho)^{2s}, t_o+6(4\varrho)^{2s}\big]\subset E_T.
\]
\end{theorem}
\noi When $s\in[\frac12,1)$, the above theorem is strictly stronger than a decisive result of Kassmann \& Weidner for the drift-free case; see \cite[Theorem~1.9]{Kass-23}. It is unclear whether such kind of result can be shown via a known adaptation of Moser's approach because of the drift term.

The next result is a dimensionless version of the previous two theorems. Note carefully that we do not assume \eqref{Eq:b:1} now, and $b$ is only required to be sufficiently integrable and does not necessarily vanish for $s\in(0,\frac12)$.
\begin{theorem}\label{Prop:WHI:3}
Let $u$ be a  local, weak super-solution to \eqref{Eq:1:1}   in $K_4 \times (0,4]$, such that $u\ge0$. 
Assume that for some $C_3>0$ it holds
$$
\|b\|_{L^{\infty}(0,4;L^{\frac{N}s}(K_2))}\le C_3.
$$
There exists a constant   
 $\eta\in(0,1)$ depending on the data $\{s, N, C_o, C_1\}$ and also on $C_3$, such that 
\begin{equation*}
	\essinf_{K_1\times(2,4]} u\ge\eta \int_{0}^1\int_{\rn}\frac{u(x,t)}{1+|x|^{N+2s}}\,\dx\dt.
\end{equation*}
\end{theorem}

\noi The above result is  new for all $s\in(0,1)$. 
Once we have these weak Harnack estimates at our disposal, together with a properly selected scaling,
deriving sub-potential lower bounds for fundamental solutions
should be the next step, cf.~\cite{BJ-06, MM-13-jfa, MM-13-aim}. We will pursue this topic in a forthcoming work.

The time-gap in weak Harnack estimates was deemed a characteristic. However, very recently, it has been understood that it is unnecessary for {\it global solutions} in the drift-free case.
That is, the following ``time-insensitive" version holds true (see~\cite[Theorem~1.1]{LW-24}):
\[
	\essinf_{K_1} u(\cdot, 1)\ge\eta \int_{\rn}\frac{u(x,1)}{1+|x|^{N+2s}}\,\dx.
\]
It would be interesting to know if such kind of ``time-insensitive" estimate still holds in the drift-diffusion case.

\subsection{Local boundedness}
Now we switch to local properties of sub-solutions.
The following local boundedness estimate generalizes \cite[Theorem~1.8]{Kass-23} in three ways.
First of all, our result handles the drift term. Second, the integral exponent is reduced from $2$ to any $\varep\in(0,2]$ for {\it sub-solutions} not just for {\it solutions}. Third, we present the estimate in a general cylinder $Q(\rho,\theta)$ instead of a parabolic cylinder $Q_\rho$. This kind of formulation, though local, is crucial in characterizing global behavior of non-negative solutions in $\rn\times(0,T]$, cf.~\cite[Chapters~V \& XI]{DB}.
It is our intension to carry out detailed calculations and to bring various points to light. Later in \S~\ref{S:discussion} we also explain that the tail term and the $L^{\varep}$-term cannot be interpolated in general, thus ruling out a conventional trick of obtaining nonlocal Harnack estimates.

\begin{theorem}\label{Prop:A:1}
 Let $u$ be a  local, weak sub-solution to \eqref{Eq:1:1}  in $E_T$. Assume that $b$ satisfies \eqref{Eq:b:1}. Then $u_+$ is locally bounded in $E_T$. Moreover, suppose $(x_o,t_o)+Q(\rho,\theta)\subset E_T$,  and let $\sig\in(0,1)$ and $\varep\in(0,2]$.  There exist a constant $\boldsymbol\gm_\varep>1$ depending only on the data $\{s, N, C_o, C_1\}$, $C_2$ and $\varep$, and a constant $q>1$ depending only on $\{s, N\}$, such that
\begin{align*}
\essup_{Q(\sig\rho,\sig \theta)} u &\le \Big(\frac{4}{1-\sig}\Big)^{N+2s} {\rm Tail} [u_+; Q(\sig \rho, \theta)]\\
&\qquad+\frac{\dsty\boldsymbol\gm_{\varep}   \Big(\frac{\rho^{2s}}{\theta}+\frac1{\sig^{N+2s}}\Big)^{\frac{2q}{\varep}}\Big(\frac{\theta}{\rho^{2s}} \Big)^{\frac1\varep}  }{(1-\sig)^{ q(N+2) -\frac{N+1}{2}}}  \Big[\frac1{(1-\sig)^{N+1}}\biint_{Q(\rho, \theta)} u_+^\varep\,\dx\dt\Big]^{\frac1\varep}.
\end{align*}
Here, we omitted the reference to $(x_o,t_o)$ in the estimate and defined
\begin{equation*}
q:=\left\{
\begin{array}{cl}
\frac{N+2s}{4s} \quad&\text{if}\quad 2s<N,\\[5pt]
\in(1,\infty) \quad &\text{if} \quad 2s\ge N.
\end{array}\right.
\end{equation*} 
\end{theorem}

An immediate yet useful, dimensionless version of the above boundedness estimate is in order.
Note again that $b$ does not have to vanish even for $s\in(0,\frac12)$.
\begin{theorem}\label{Prop:A:2}
Let $u$ be a  local, weak sub-solution to \eqref{Eq:1:1}   in $K_4 \times (0,4]$. Assume that for some $C_3>0$ it holds
$$
\|b\|_{L^{\infty}(0,4;L^{\frac{N}s}(K_2))}\le C_3.
$$
There exists a constant   
 $\boldsymbol\gm>1$ depending on the data $\{s, N, C_o, C_1\}$ and also on $C_3$, such that 
\begin{equation*}
\essup_{K_{1}\times(3, 4]}u\le \boldsymbol\gm \int_{0}^4\int_{\rn}\frac{u_+(x,t)}{1+|x|^{N+2s}}\,\dx\dt.
\end{equation*}
\end{theorem} 

\noi We stress that the above two results hold for {\it sub-solutions} not just for {\it solutions}. In \cite{Kass-23} for the drift-free case and in \cite{NNSW} for the drift case, an $L^2-L^{\infty}$ estimate is first established for {\it sub-solutions}, and then an $L^1-L^{\infty}$ estimate for {\it solutions} is derived from it using a covering argument. 

Very recently, we were able to improve this kind of estimate and obtain a ``time-insensitive" version of it for {\it global solutions} in the case $b=0$. Namely, we have (see~\cite[Theorem~1.1]{LW-24})
\[
\essup_{K_{1}}u(\cdot, 1)\le \boldsymbol\gm \int_{\rn}\frac{u_+(x,1)}{1+|x|^{N+2s}}\,\dx.
\]
A natural question would be if it still holds for the drift-diffusion case.

\subsection{Harnack's estimates}
Harnack's inequality for non-negative solutions now follows from a combination of Theorems~\ref{Prop:WHI:2-}, \ref{Prop:WHI:2} and \ref{Prop:A:1}.
\begin{theorem}\label{Prop:HI:1}
Let $u$ be a  local, weak solution to \eqref{Eq:1:1}   in $E_T$, such that $u\ge0$. Assume that $b$ satisfies \eqref{Eq:b:1}.
There exists a constant   
 $\boldsymbol\gm>1$ depending on the data $\{s, N, C_o, C_1\}$ and also on $C_2$, such that 
\begin{equation*}
	\essup_{K_{\rho}(x_o)\times( t_o-\varrho^{2s},
	t_o)} u\le \boldsymbol\gm  \essinf_{K_{\rho}(x_o)\times( t_o+3\varrho^{2s},
	t_o+4\varrho^{2s})} u
\end{equation*}
provided that  
\[
K_{4\rho}(x_o)\times\big(t_o-4\varrho^{2s}, t_o+4\varrho^{2s}\big]\subset E_T.
\]
\end{theorem}
\noi When $s\in[\frac12,1)$, the above theorem is strictly stronger than the result of Kassmann \& Weidner for the drift-free case; see \cite[Theorem~1.1]{Kass-23}.

A dimensionless version of the Harnack inequality can be derived from Theorems~\ref{Prop:WHI:3} and \ref{Prop:A:2}.
Let us emphasize that $b$ only needs to be sufficiently integrable and does not necessarily vanish even if $s\in(0,\frac12)$.
\begin{theorem}\label{Prop:HI:2}
Let $u$ be a  local, weak solution to \eqref{Eq:1:1}   in $K_4 \times (-4,4]$, such that $u\ge0$. 
Assume that for some $C_3>0$ it holds
$$
\|b\|_{L^{\infty}(-4,4;L^{\frac{N}s}(K_2))}\le C_3.
$$
There exists a constant   
 $\boldsymbol\gm>1$ depending on the data $\{s, N, C_o, C_1\}$ and also on $C_3$, such that 
\begin{equation*}
	 \boldsymbol\gm^{-1} \essup_{K_{1}\times( -1,0)} u\le \int_{-2}^0\int_{\rn}\frac{u(x,t)}{1+|x|^{N+2s}}\,\dx\dt\le  \boldsymbol\gm  \essinf_{K_1\times( 3,4)} u.
\end{equation*}
\end{theorem}
\noi Note that the above Harnack's inequality is  new for all $s\in(0,1)$ due to the presence of the drift term. It also offers some leeway to investigate the less understood range $s\in(0,\frac12)$ in the future.
With a properly chosen scaling, we expect that it will tell more precise local and global behavior of solutions.

Here, Harnack estimates require a time-gap because {\it local solutions} are considered. Recently in \cite{LW-24}, it has been understood that, for the drift-free case, such a time-gap is unnecessary for {\it global solutions} and {\it elliptic-type} Harnack estimates actually hold. A natural question now is whether this continues to hold for drift-diffusion equations. A positive answer seems to hint immediate regularization of solutions.
\subsection{H\"older regularity}
A DeGiorgi-type approach was successfully implemented in the study of H\"older regularity for nonlocal drift-diffusion equations in \cite{Caff-Vass-10}; see also \cite{Caff-Vass-11}.
Since then it has been improved in different aspects, cf. \cite{Const-Wu, DS-18, NSL, Silv-10} for instance. For a viscosity approach, see~\cite{Silv-11, Silv-12, Silv-12p} for instance.
To obtain H\"older regularity for nonlocal, nonlinear diffusion equations, a DeGiorgi-type scheme was devised in \cite{Liao-cvpd-24} which we also follow here. Such a scheme has already been implemented in \cite{NNSW} to improve H\"older estimates for the critical case $s=\frac12$ of \eqref{Eq:1:1}; see also \cite{Byun-ampa, Liao-mod-24} for some refinements of the scheme in the case of nonlocal, nonlinear diffusion equations.
Our use of DeGiorgi's method follows closely the interpretation of Lady\v{z}enskaja \&  Ural'tseva (see \cite{LSU} for instance) and is also influenced by DiBenedetto. 

Here, we present an improved H\"older estimate for the nonlocal drift-diffusion equation \eqref{Eq:1:1}
in the so-called sub-critical range of $s$.
\begin{theorem} \label{Thm:holder:1}
Let $s\in(\frac12,1)$ and $u$ be a locally bounded, local, weak solution to \eqref{Eq:1:1} in $E_T$.
Assume that $b$ satisfies \eqref{Eq:b:1} and that for some $\varep>0$,
	\begin{equation}\label{Eq:global-thm}
	\int_{\rn}\frac{|u(x,\cdot)|}{1+|x|^{N+2s}}\,\dx\in L^{1+\varep}_{\loc}(0,T].
	\end{equation}
	Then $u$ is locally H\"older continuous in $E_T$. More precisely,
	there exist a constant $\boldsymbol\gm>1$ depending on the data $\{s, N, C_o, C_1\}$ and a constant $\be\in(0,1)$ depending on the data and $\varep$, such that for any $0<r<R<\widetilde{R}$, there holds
	\begin{equation*}
	\essosc_{(x_o,t_o)+Q_r}u \le  \boldsymbol\gm\boldsymbol\,\om \Big(\frac{r}{R}\Big)^{\be},
	\end{equation*}
provided the  cylinders  $(x_o,t_o)+Q_{R} \subset (x_o,t_o)+Q_{\widetilde{R}}$ 
are included in $E_T$, 
where
$$\boldsymbol\om=2\essup_{(x_o,t_o)+Q_{\widetilde{R}}}|u| +\bigg(\bint^{t_o}_{t_o- \widetilde{R}^{2s}}  \Big( \widetilde{R}^{2s} \int_{\rn\setminus K_{\widetilde{R}}(x_o)} \frac{|u(x,t)| }{|x-x_o|^{N+2s}}\,\dx\Big)^{ 1+\varep} \dt\bigg)^{\frac1{1+\varep}}.$$
In particular, if $b\in L^{\infty}(0,T;L_{\loc}^{\frac{N}{2s-1}}(E_T))$, then the same oscillation estimate holds in dependence of the corresponding norm of $b$.
\end{theorem}

\noi This improves \cite[Theorems~1.1 \& 1.4]{NNSW}, where $b(\cdot,t)$ is required to be bounded or in BMO.
Such kind of result seems not achievable via a known adaptation of Moser's approach; see \cite[Remark~2.1]{NNSW} in this connection.

As remarked after \cite[Theorem~1.1]{NNSW}, this kind of result is new in contrast to previous H\"older estimates for nonlocal drift-diffusion equations because of two reasons. First, the estimate is truly local in nature in the sense that it does not require $u$ to solve an equation globally and dispenses with any kind of harmonic extension. Second, our  H\"older regularity only requires local boundedness, and solutions can grow at infinity as long as the tail term is finite. Naturally, one wonders if local continuity of $u$ continues to hold for the borderline case $\varep=0$ in \eqref{Eq:global-thm}; see \cite{Liao-mod-24} for the drift-free case.

The above H\"older regularity could have been stated uniformly for all $s\in(0,1)$ under the condition \eqref{Eq:b:1} of $b$.
However, for $s\in(0,\frac12]$ one can properly choose a different scaling and carry out more refined analysis under the following weaker condition:  
\begin{equation*}
\essup_{t\in[0,T]} \Big(\bint_{K_{R}(x_o)} \big|b(x,t)-[b(\cdot, t)]_{K_R(x_o)}\big|^{\frac{N}s}\,\dx\Big)^{\frac{s}N}\le C_2 R^{1-2s}
\end{equation*}
holds true for all $K_R(x_o)\subset E$.
When $s=\frac12$, such refined analysis has been performed in \cite{NNSW} following a particular scaling from \cite{Caff-Vass-10}. We believe that with a proper scaling, the same program can also be applied for $s\in(0,\frac12)$ and improve known H\"older regularity in this range as well, cf.~\cite{Const-Wu, NSL}. Needless to say, the regularity of solutions to the surface quasi-geostrophic equation with $s$ in this range is still a major open problem.

\subsection{Technical innovation}
 
To establish weak Harnack estimates, an important feature of our approach is that, equation \eqref{Eq:1:1} is used only once to derive energy estimates. All our results are sole consequences of the energy estimates, and the equation does not play a role any more. Thus, the results are quite robust. Such point of view is in line with our previous works~\cite{Liao-cvpd-24, Liao-mod-24}, whose theme is the local continuity of solutions.  
This feature is distinct from a Moser-type scheme, which keeps using the equation with different testing functions and usually hinges upon certain crossover estimates based on the BMO estimate of $\log u$. As such, our method is more flexible in dealing with a drift-term, and also holds the promise of treating more general objects like quasi-minima from Calculus of Variations.

Even though  a DeGiorgi-type scheme was first applied in \cite{Caff-Vass-10} to nonlocal drift-diffusion equations to obtain H\"older regularity, it is well known that Harnack estimates are stronger than H\"older regularity. To establish Harnack's inequality, no matter for local or nonlocal operators, under the DeGiorgi-Nash-Moser theory, most literatures follow Moser's approach. If one instead follows DeGiorgi's approach, then certain covering argument of Krylov-Safanov is usually needed; see~\cite{DB-Trud}. Therefore, to achieve nonlocal weak Harnack estimates for super-solutions, more refined local and nonlocal estimates than those in \cite{Caff-Vass-10, Caff-Vass-11} are asked for. In fact, our Theorem~\ref{Prop:WHI:2-} includes a first,  DeGiorgi-type proof of a nonlocal, parabolic weak Harnack estimate obtained by Felsinger \& Kassmann \cite{Kass-13}.  Amongst the technical tools used, we would like to highlight a DeGiorgi-type lemma \ref{Lm:DG:initial:1}, which is new even for local operators, and a new, nonlocal, measure shrinking mechanism (Lemmata~\ref{Lm:3:2} \& \ref{Lm:3:2+}) that allow us, differently from the local operator theory, to avoid a Krylov-Safanov-type covering argument.  Also because of them, we are able to handle the difficulty of certain critical scaling property brought by the tail.

To prove a refined $L^{\infty}$-estimate for sub-solutions, we used a time-dependent truncation level to balance the nonlocal tail in the DeGiorgi iteration. This trick is taken from  \cite{Kass-23}. However, we point out that on a different occasion, a time-dependent truncation level appeared earlier in \cite[Proposition~2.1]{DGV-07} in order to obtain refined estimates from DeGiorgi's iteration.

We expect our approach is extendable and has the potential to be applied to study nonlocal, nonlinear parabolic equations, as evidenced in the local operator theory, cf.~\cite{DB,DBG-pde, DBGV-mono}; indeed, how to interpret the existing local, nonlinear parabolic theory via Moser's approach is still an interesting open problem.

Finally we point out that, our measure theoretical approach captures some peculiar characters of the nonlocal equations that were overlooked in the literature even in the elliptic case, and as a result, it also readily yields novel findings in this case. 
For instance, 
we will show in \S~\ref{S:elliptic} that, if $u$ is globally non-negative and belongs to a certain fractional DeGiorgi class ${\rm DG}_-^{s,p}$, which is modeled on super-solutions to the elliptic fractional $p$-Laplacian, then there exists $\eta\in(0,1)$, such that
\begin{equation}\label{Eq:intro:3}
\essinf_{K_{1} } u  \ge \eta\Big( \int_{\rn }\frac{u^{p-1}(x)}{1+|x|^{N+sp}}\,\dx \Big)^{\frac1{p-1}}.
\end{equation}
The novelty is two-fold: The positivity contribution from the long-range behavior of $u$ is discovered and the (local) integral exponent is raised to $p-1$ matching perfectly the integral exponent of the nonlocal tail. The significance of such improvement lies in potential applications to the Wiener-type criterion for quasi-minima \cite{DBG-Wiener}; we refer the reader to \S~\ref{N-S-p} for more discussion.

\subsection{Definitions and notation}\label{S:notion}

\subsubsection{Notion of solution}
\begin{definition}\label{Def:solution}
A measurable function $u:\,\rn\times(0,T]\to\rr$ satisfying
\begin{equation*} 
	u\in C_{\loc}\big(0,T;L^2_{\loc}(E)\big)\cap L^2_{\loc}\big(0,T; W^{1,2}_{\loc}(E)\big)
\end{equation*}
is a local, weak sub(super)-solution to \eqref{Eq:1:1}, if for every bounded open set $\Om\subset E$ and every sub-interval
$[t_1,t_2]\subset (0,T]$, we have
\begin{equation}\label{Eq:global-int}
\int_{t_1}^{t_2}\int_{\rn}\frac{|u(x,t)|}{1+|x|^{N+2s}}\,\dx\dt<\infty
\end{equation}
and
\begin{equation*} 
\begin{aligned}
	\int_{\Om} & u\vp \,\dx\bigg|_{t_1}^{t_2}
	-\int_{t_1}^{t_2}\int_{\Om}  u\big(\pl_t\vp+b\cdot\nabla\vp \big)\, \dx\dt
	+\int_{t_1}^{t_2} \mathscr{E}\big(u(\cdot, t), \vp(\cdot, t)\big)\,\dt
	\le(\ge)0
\end{aligned}
\end{equation*}
where
\[
	 \mathscr{E}:=\int_{\rn}\int_{\rn}a(x,y,t)  \frac{ \big(u(x,t) - u(y,t)\big)\big(\vp(x,t) - \vp(y,t)\big)}{2|x-y|^{N+2s}}\,\dy\dx
\]
for all non-negative testing functions
\begin{equation*} 
\vp \in W^{1,2}_{\loc}\big(0,T;L^2(\Om)\big)\cap L^2_{\loc}\big(0,T;W_o^{1,2}(\Om)%
\big).
\end{equation*}

A function $u$ that is both a local weak sub-solution and a local weak super-solution
to \eqref{Eq:1:1}  is a local weak solution.

\end{definition}

\begin{remark}\upshape
For sub/super-solutions, it is possible to weaken $C_{\loc}$ in time to $L^{\infty}_{\loc}$, whereas solutions are automatically in the former space because of their weak formulation.
\end{remark}

\subsubsection{Notation}
Let $K_\rho(x_o)$ be the ball of radius $\rho$ and center $x_o$ in $\rn$.
We use  the symbols 
\begin{equation*}
\left\{
\begin{aligned}
(x_o,t_o)+Q(R,S)&:=K_R(x_o)\times (t_o-S,t_o),\\[5pt]
(x_o,t_o)+Q_\rho(\theta)&:=K_{\rho}(x_o)\times(t_o-\theta\rho^{2s},t_o),
\end{aligned}\right.
\end{equation*} 
to denote (backward) cylinders.
The vertex $(x_o,t_o)$ is often omitted from the cylinder in \eqref{Eq:tail} for simplicity. 
If $\theta=1$, it will also be omitted. When the context is unambiguous, we will apply these conventions.

A central object of our study is the tail:
\begin{equation}\label{Eq:tail}
{\rm Tail}[u; Q(R,S)]:= \int_{t_o-S}^{t_o} \int_{\rn\setminus K_R(x_o)}\frac{|u(x,t)|}{|x-x_o|^{N+2s}}\,\dx \dt.
\end{equation}
Its finiteness is ensured by \eqref{Eq:global-int}.

Throughout this paper,  we use $\boldsymbol\gm$ as a generic positive constant in various estimates that can be determined by the data from the problem under consideration.

\subsection{Organization}
Section~\ref{S:energy} is devoted to the energy estimates. Section~\ref{S:prelim} collects main modules needed in the proof of weak Harnack estimates, which are presented in Section~\ref{S:Harnack}. In Section~\ref{S:bd}, we derive a local boundedness estimate for sub-solutions, whereas H\"older regularity is proven in Section~\ref{S:Holder}. Finally, in Section~\ref{S:elliptic} an application of our approach in an elliptic setting is demonstrated.

\

\noi{\bf Acknowledgement.} This work was supported by the FWF-project P36272-N ``On the Stefan type problems."

\section{Energy estimates}\label{S:energy}

The energy estimates in \cite[Proposition~2.1]{Liao-cvpd-24} \& \cite[Proposition~2.1]{Liao-mod-24} are sufficient for studying the local continuity of solutions. However, to deal with Harnack estimates we need to enhance these estimates, i.e. an integral over the whole $\rn$ appears on the left-hand side.  
In fact, the following energy estimates practically encode all the information needed to show H\"older regularity and weak Harnack estimates. 
\begin{proposition}\label{Cor:1}
	Let $u$ be a  local weak sub(super)-solution to \eqref{Eq:1:1} in $E_T$.
	There exists a constant $\boldsymbol \gm (C_o,C_1)>0$, such that
 	for all cylinders $Q(r,\tau)\subset Q(R,S) \subset E_T$ with their vertex at $(x_o,t_o)$,
 	and every $k\in\rr$, 
	there holds
\begin{align}\nonumber
\essup_{t_o-\tau<t<t_o}&\int_{K_r} w^2_{\pm}(x,t)\,\dx+
	\int_{t_o-\tau}^{t_o}\int_{K_r}\int_{K_r}  \frac{|w_\pm(x,t) - w_\pm(y,t)|^2}{|x-y|^{N+2s}}\,\dx\dy\dt\\ \nonumber
	&\le
	\frac{\boldsymbol\gm R^{2(1-s)}}{(R-r)^2}\iint_{Q(R,S)}  w^2_{\pm}(x,t) \,\dx\dt+ \frac{\boldsymbol\gm}{S-\tau}\iint_{Q(R,S)} w_{\pm}^2(x,t)\,\dx\dt \\ \label{Eq:cor:1}
	&\quad+\frac{\boldsymbol\gm R^{N}}{(R-r)^{N+2s}}\iint_{Q(R,S)} w_{\pm}(x,t)\,\dx\dt  
	\int_{\rn\setminus K_{\frac{r+R}2}}\frac{w_{\pm}(y,t)}{|y-x_o|^{N+2s}}\,\dy
\end{align}
and 
\begin{align}\nonumber
	 \essup_{t_o-S<t<t_o} &\int_{K_r} w^2_{\pm}(x,t)\,\dx + \int_{t_o-S}^{t_o}\int_{K_r}\int_{K_r}  \frac{|w_\pm(x,t) - w_\pm(y,t)|^2}{|x-y|^{N+2s}}\,\dx\dy\dt\\ \nonumber
	 &\quad+\iint_{Q(R,S)} w_{\pm}(x,t)\,\dx\dt \bigg(\int_{ \rn}  \frac{w_\mp(y,t)}{|x-y|^{N+2s}}\,\dy\bigg)  \\ \nonumber
	 &\le \int_{K_R} w^2_\pm(x,t_o-S)\,\dx
	+\frac{\boldsymbol\gm R^{2(1-s)}}{(R-r)^2} \iint_{Q(R,S)}  w^2_{\pm}(x,t) \,\dx\dt\\ \label{Eq:cor:2}
	&\quad +\frac{\boldsymbol\gm R^{N}}{(R-r)^{N+2s}}\iint_{Q(R,S)} w_{\pm}(x,t)\, \dx\dt
	\int_{\rn\setminus K_{\frac{r+R}2}}\frac{w_{\pm}(y,t)}{|y-x_o|^{N+2s}}\,\dy.
\end{align}
Here, we have denoted $w=u-k$ for simplicity.
\end{proposition}
 \begin{proof}
 We only deal with sub-solutions. Let $\z\in C^1_o(Q(R,S))$ be a non-negative cutoff function. 
 Use $w_+\z^2$ a testing function against \eqref{Eq:1:1}.
The treatment of the time part is quite standard. 
Let us focus on the diffusion part. Omitting the time variable for simplicity, we observe, due to $ w_- w_+=0$, that
\begin{align*} 
&(u(y)-u(x))[(u(y)-k)_{+}\zeta^2(y)-(u(x)-k)_{+}\zeta^2(x)]\nonumber\\
&\qquad=[\zeta(y)w_+(y)-\zeta(x)w_+(x)]^2+w_-(y)w_+(x)\zeta^2(x)+w_-(x)w_+(y)\zeta^2(y)\nonumber\\
&\qquad\qquad-w_+(y)w_+(x)(\zeta(x)-\zeta(y))^2.
\end{align*}
Omitting the time integration, this readily yields the estimate of the diffusion term, that is,
\begin{align*}
      \iint_{\rn\times\rn} & (u(y)-u(x))[(u(y)-k)_{+}\zeta^2(y)-(u(x)-k)_{+}\zeta^2(x)]\,\d\mu\\
      &\ge \iint_{\rn\times\rn}[\zeta(y)w_+(y)-\zeta(x)w_+(x)]^2\,\d\mu\\
      &\quad + \iint_{\rn\times K_R}w_-(y)w_+(x)\zeta^2(x)\, \d\mu \\
      &\quad - \iint_{\rn\times\rn} w_+(y)w_+(x)(\zeta(x)-\zeta(y))^2\, \d\mu
\end{align*}
where we set $\d\mu=K(x,y,\cdot)\dy\dx$.
The last term can be rewritten as
\begin{align*}
      \iint_{\rn\times\rn} & w_+(y)w_+(x)(\zeta(x)-\zeta(y))^2\, \d\mu\\
      &= \iint_{K_R\times K_R}  w_+(y)w_+(x)(\zeta(x)-\zeta(y))^2\, \d\mu \\
      &\quad +2 \iint_{K_R\times (\rn\setminus K_R)} \z^2 w_{+}(x)w_{+}(y)\,\d\mu.
\end{align*}
 The drift term can be treated as follows. Formally, we have
 \begin{align*}
      -\iint_{Q(R,S)} \z^2  b\cdot \nabla w w_{+}\,\dx\dt &= -\frac12 \iint_{Q(R,S)} \z^2  b\cdot \nabla w^2_{+}\,\dx\dt\\
      &=\frac12\iint_{Q(R,S)} \dvg (\z^2  b)   w^2_{+}\,\dx\dt\\
      &=\iint_{Q(R,S)} \z\nabla\z \cdot b   w^2_{+}\,\dx\dt      
 \end{align*}
 where we used  $\dvg b(\cdot, t)=0$ in the last line. We continue to estimate the last term by Young's inequality and then by Sobolev-type embedding and H\"older's inequality as
 \begin{align*}
      \iint_{Q(R,S)} \z |\nabla\z| |b|   w^2_{+}\,\dx\dt &\le \varep \int_{t_o-S}^{t_o}\Big(\int_{K_R}(\z w_+)^{\frac{2N}{N-2s}}\,\dx\Big)^{\frac{N-2s}{N}}\dt\\
      &\quad+\boldsymbol\gm(\varep)\int_{t_o-S}^{t_o}\Big(\int_{K_R}(|\nabla\z| |b| w_+)^{\frac{2N}{N+2s}}\,\dx\Big)^{\frac{N+2s}{N}}\dt\\
      &\le \varep C_{\mathrm{Sob}}\int_{t_o-S}^{t_o} \iint_{\rn\times\rn}[\zeta(y)w_+(y)-\zeta(x)w_+(x)]^2\,\d\mu \dt\\
      &\quad +\boldsymbol\gm(\varep)\int_{t_o-S}^{t_o}\Big(\int_{K_R} |b|^{\frac{N}s}\,\dx\Big)^{\frac{2s}N}\Big(\int_{K_R}|\nabla\z|^2 w_+^2\,\dx\Big)\,\dt.
 \end{align*}
 The first term on the right-hand side can be absorbed to the left-hand side if we select $\varep=1/(2C_{\mathrm{Sob}})$, whereas the integral of $b$ in the second term is estimated by \eqref{Eq:b:1}. The desired energy estimates follow from the above estimates and properly choosing $\z$.
 \end{proof}
The following energy estimate will be used to derive local boundedness estimates.
We use a time-dependent truncation method from \cite{Kass-23}.  
\begin{proposition}\label{Prop:2:1}
	Let $u$ be a  local weak sub(super)-solution to \eqref{Eq:1:1}  in $E_T$, and
	let $k(\cdot)$ be  absolutely continuous  in $(0,T)$.
	There exists a constant $\boldsymbol \gm (C_o,C_1)>0$, such that
 	for all cylinders $Q(R,S) \subset E_T$,
 	and every non-negative, piecewise smooth cutoff function
 	$\z(\cdot,t)$ compactly supported in $ K_{R} $ for all  $t\in (t_o-S,t_o)$,  there holds
\begin{align*}
	&\int_{t_o-S}^{t_o}\int_{K_R}\int_{K_R}  \frac{|\z w_\pm(x,t) - \z w_\pm(y,t)|^2}{|x-y|^{N+2s}}\,\dx\dy\dt\\
	&\qquad+\iint_{Q(R,S)} \z^2 w_{\pm}(x,t)\,\dx\dt \bigg(\int_{ \rn}  \frac{w_\mp(y,t)}{|x-y|^{N+2s}}\,\dy\bigg)
	 +\int_{K_R}\z^2 w^2_{\pm}(x,t)\,\dx\bigg|_{t_o-S}^{t_o}\\
	&\quad\le
	\boldsymbol\gm\int_{t_o-S}^{t_o}\int_{K_R}\int_{K_R}\max\big\{w^2_{\pm}(x,t), w^2_{\pm}(y,t)\big\} \frac{|\z(x,t) - \z(y,t)|^2}{|x-y|^{N+2s}}\,\dx\dy\dt\\
	&\qquad+\boldsymbol\gm\int_{t_o-S}^{t_o}\int_{K_R}\int_{\rn\setminus K_R} \z^2 w_{\pm}(x,t)\frac{ w_{\pm}(y,t)}{|x-y|^{N+2s}}\,\dy\dx\dt 
	\\
	&\qquad \mp 2\iint_{Q(R,S)}  k'(t) \z^2 w_{\pm}(x,t)\,\dx\dt + \iint_{Q(R,S)} |\pl_t\z^2 |w_{\pm}^2(x,t)\,\dx\dt. 
\end{align*}
Here, we have denoted $w(x,t)=u(x,t)-k(t)$ and omitted the vertex $(x_o,t_o)$ for simplicity.
\end{proposition}
\begin{proof}
Take the case of sub-solution for instance. Using $\vp=w_+\z^2$ as a testing function in the weak formulation, the integrals resulting from the fractional diffusion term and the drift term are treated as in Proposition~\ref{Cor:1}. Regarding the part with the time derivative, formally we can write
\[
 \pl_t u \cdot (u-k)_+\z^2 =\tfrac12  \z^2\pl_t (u-k)_+^2  + \z^2 k^{\prime}(u-k)_+ .
\]
Integrating this equality in $Q(R,S)$ readily yields the conclusion. A rigorous treatment of the time derivative can be adapted from 
\cite[Appendix~B]{Liao-cvpd-24}.
\end{proof}

\section{Preparatory materials} \label{S:prelim}

In this section, we collect main modules for the proof of Harnack estimates. 
Some of them have already appeared in \cite[Section~3]{Liao-mod-24}. However, Lemma~\ref{Lm:DG:initial:1}, Lemma~\ref{Lm:3:2} and Lemma~\ref{Lm:3:2+} are brand-new.

\subsection{DeGiorgi-type lemmata}
Based on the energy estimate \eqref{Eq:cor:1} the following DeGiorgi-type lemma appears in \cite[Lemma~3.1]{Liao-mod-24}.
\begin{lemma}\label{Lm:DG:1}
 Let $u$ be a  local weak super-solution to \eqref{Eq:1:1}  in $E_T$, such that $u\ge0$.
 Let $k>0$ and $\theta\in(0,1)$ be parameters.
There exist  constants $\widetilde{\boldsymbol\gm}>1$ depending only on 
 the  data $\{s, N, C_o, C_1\}$, and   $\nu\in(0,1)$ depending on the data $\{s, N, C_o, C_1\}$ and $\theta$, such that if
\begin{equation*}
	\big|\big\{
	u\le k\big\}\cap  (x_o,t_o)+Q_{\varrho}(\theta)\big|
	\le
	\nu|Q_{\varrho}(\theta)|,
\end{equation*}
then 
\begin{equation*}
	u\ge \tfrac14k
	\quad
	\mbox{a.e.~in $ (x_o,t_o)+Q_{\frac{1}2\varrho}(\theta)$,}
\end{equation*}
provided that $(x_o,t_o)+ Q_{2\varrho}(\theta) \subset E_T$.
Moreover,  $\nu= \bar\nu \theta^q$ for some $\bar\nu\in(0,1)$ depending only on the data $\{s, N, C_o, C_1\}$, and $ q>1$ depending only on $N$.
\end{lemma}

The previous lemma says that, if a critical measure density is reached in a cylinder, then pointwise information can be extracted in a smaller cylinder attached to the same vertex. For this reason, it is also referred to as the {\it critical mass/density lemma} in the literature.
 The next lemma, like the previous one, features a {\it from-measure-to-pointwise} estimate.  It asserts that, if a critical measure density is reached in a spatial ball at some time, then pointwise information can be claimed at later times. The crucial difference from the previous one lies in how the measure theoretical information is assigned.
As such, it constitutes the first new element of this section.
\begin{lemma}\label{Lm:DG:initial:1}
Let $u$ be a   local weak super-solution to \eqref{Eq:1:1}   in $E_T$, such that $u\ge0$. Let $k>0$ be a parameter. 
There exist constants $\dl,\,\nu_o\in(0,1)$ and $\widetilde{\boldsymbol\gm}>1$ depending only on the data $\{s, N, C_o, C_1\}$, such that if
\[
\big| \big\{u(\cdot, t_o)\le k\big\}\cap  K_{2\rho}(x_o)\big|<\nu_o |K_{2\rho}|,
\]
then 
\[
u\ge \tfrac1{8} k\quad\text{ a.e. in }K_{\frac12\rho}(x_o)\times\big(t_o+\tfrac34\dl \rho^{2s}, t_o+\dl \rho^{2s}\big],
\]
provided that $K_{2\rho}(x_o)\times(t_o, t_o+\dl\rho^{2s}]\subset E_T$.
\end{lemma}
\begin{proof}
 Upon a translation, we may assume $(x_o,t_o)=(0,0)$ and
define for $n\in \nn_0$,
\begin{align*}
	\left\{
	\begin{array}{c}
\dsty k_n= 2^{-1}k +2^{-(n+1)} k,
\\[5pt]
\dsty \rho_n=\rho+2^{-n}\rho, \quad\widetilde{\rho}_n=\tfrac12(\rho_n+\rho_{n+1}),\\[5pt]
\dsty \widehat{\rho}_n=\tfrac34\rho_n+ \tfrac14\rho_{n+1},\quad \overline{\rho}_n=\tfrac14\rho_n+ \tfrac34\rho_{n+1},\\[5pt]
\dsty K_n=K_{\rho_n}, \quad \widetilde{K}_n=K_{\widetilde{\rho}_n},\quad \widehat{K}_n=K_{\widehat\rho_n},\quad \overline{K}_n=K_{\overline{\rho}_n},\\[5pt]
\dsty Q_n=K_n\times(0,\dl\rho^{2s}],\quad\widetilde{Q}_n=\widetilde{K}_n\times(0,\dl\rho^{2s}],\\[5pt]
\widehat{Q}_n=\widehat{K}_n\times(0,\dl\rho^{2s}],\quad\overline{Q}_n=\overline{K}_n\times(0,\dl\rho^{2s}].
\end{array}
	\right.
\end{align*}
Observe that $Q_{n+1}\subset\overline{Q}_n\subset\widetilde{Q}_n\subset\widehat{Q}_n\subset Q_n$ and their height is fixed to be $\dl\rho^{2s}$, where $\dl$ is to be chosen.
The energy estimate \eqref{Eq:cor:2} of Proposition~\ref{Cor:1} is used in $\widetilde{Q}_n\subset\widehat{Q}_n$ with 
 $$
w_-(x,t):=\big(u(x,t) -k_n\big)_-. 
$$
This gives 
\begin{align*}
	\max\bigg\{ &\essup_{t\in[0,\dl\rho^{2s}]}\int_{\widetilde{K}_n} w^2_{-}(x,t)\,\dx,\, \int^{\dl\rho^{2s}}_{0}\int_{\widetilde{K}_n}\int_{\widetilde{K}_n}  \frac{|w_-(x,t) - w_-(y,t)|^2}{|x-y|^{N+2s}}\,\dx\dy\dt \bigg\}\\
	&\quad\le  \int_{K_n}  w_{-}^2(x,0)\,\dx+
	\frac{\boldsymbol\gm 4^n}{\rho^{2s}}\int^{\dl\rho^{2s}}_{0}\int_{K_n} w^2_{-}(x,t) \,\dx\dt\\
	&\qquad+\boldsymbol\gm\iint_{Q_n}  w_{-}(x,t)\,\dx\dt \bigg(  \int_{\rn\setminus \widehat{K}_n}\frac{ w_{-}(y,t)}{|y|^{N+2s}}\,\dy\bigg)
\end{align*}
Proceeding as in \cite[Lemma~3.1]{Liao-mod-24}, the last two terms on the right-hand are estimated by
\[
 \boldsymbol\gm 2^{(N+4)n}\frac{k^2}{\rho^{2s}}|A_n|,
\]
where
\[
A_n:=\big\{u<k_n\big\}\cap Q_n.
\]
Note that $\dl$ does not appear in the above estimate as $\z$ is independent of the time variable.

To estimate the first term of the energy estimate we observe that, since $\rho_n\in(\rho, 2\rho]$ and $k_n\le k$, we have
\begin{align*}
\big|\big\{u(\cdot, 0)<k_n\big\}\cap K_n\big| &\le \big|\big\{u(\cdot, 0)<k\big\}\cap K_{2\rho}\big|\\
&\le \nu_o |K_{2\rho}| \le 2^N\nu_o |K_n| =\frac{2^N\nu_o}{\dl\rho^{2s}}|Q_n|,
\end{align*}
where $\nu_o$ is still to be selected.
Then, it is estimated by
\[
\int_{K_n}  w_{-}^2(x,0)\,\dx \le k^2 \big|\big\{u(\cdot, 0)<k_n\big\}\cap K_n\big|\le 2^N\frac{k^2}{\dl\rho^{2s}}\nu_o|Q_n|.
\]
Therefore, the energy estimate becomes
\begin{equation}\label{Eq:DG-energy:1}
\begin{aligned}
	\essup_{t\in[0,\dl\rho^{2s}]}&\int_{\widetilde{K}_n} w^2_{-}(x,t)\,\dx+ \int^{\dl\rho^{2s}}_{0}\int_{\widetilde{K}_n}\int_{\widetilde{K}_n}  \frac{|w_-(x,t) - w_-(y,t)|^2}{|x-y|^{N+2s}}\,\dx\dy\dt\\
	&\quad\le
	 \boldsymbol\gm 2^{(N+4)n}\frac{k^2}{\rho^{2s}}|A_n|+2^N\frac{k^2}{\dl\rho^{2s}} \nu_o |Q_n|.
\end{aligned}
\end{equation}

To proceed, we consider two cases. In the {\it first case}, suppose
\begin{equation}\label{Eq:DG-alt:1}
|A_n|=\big|\big\{u <k_n\big\}\cap Q_n\big|\le \frac{\nu_o}{\dl}|Q_n|\quad\text{for some}\> n\in  \nn_0.
\end{equation}
Noticing that $\rho_n\in(\rho, 2\rho]$ and $k_n\ge\frac12 k$, 
then \eqref{Eq:DG-alt:1} yields that
\begin{equation}\label{Eq:meas-alt:1}
\begin{aligned}
 \big|\big\{u<\tfrac12 k\big\}\cap (0,\dl\rho^{2s})+Q_{\rho}(\dl)\big| &\le\big|\big\{u<k_n\big\}\cap Q_n\big| \\
&\le\frac{\nu_o}{\dl}|Q_n|\le 2^N\frac{\nu_o}{\dl} |Q_{\rho}(\dl)|.
\end{aligned}
\end{equation}
Recall our notation $(0,\dl\rho^{2s})+Q_{\rho}(\dl)=K_{\rho}\times(0,\dl\rho^{2s}]$ in \eqref{Eq:meas-alt:1}.
Assuming that $\dl$ has been fixed for the moment, let us choose $\nu_o$ to be so small that 
\begin{equation}\label{Eq:nu:0}
\nu_o\le 2^{-N}\bar{\nu} \dl^{q+1}
\end{equation}
where $\bar{\nu}$ and $q$ depend only on the data $\{s, N, C_o, C_1\}$, as determined in Lemma~\ref{Lm:DG:1}. Then   \eqref{Eq:meas-alt:1} and \eqref{Eq:nu:0} allow us to apply Lemma~\ref{Lm:DG:1} and conclude that
\begin{equation}\label{Eq:DG:pt-est:1}
u\ge\tfrac1{8}k\quad\text{ a.e. in }(0,\dl\rho^{2s})+Q_{\frac12\rho}(\dl).
\end{equation}
This step does not require any quantitative knowledge of $\dl$, which is still to be specified.

In the {\it second case}, suppose
\begin{equation}\label{Eq:DG-alt:2}
|A_n|=\big|\big\{u <k_n\big\}\cap Q_n\big|> \frac{\nu_o}{\dl}|Q_n|\quad\text{for any}\> n\in  \nn_0.
\end{equation}
Then, the energy estimate \eqref{Eq:DG-energy:1} gives that for any $n\in  \nn_0$,
\begin{equation*}
\begin{aligned}
	\essup_{t\in[0,\dl\rho^{2s}]}&\int_{\widetilde{K}_n} w^2_{-}(x,t)\,\dx+ \int^{\dl\rho^{2s}}_{0}\int_{\widetilde{K}_n}\int_{\widetilde{K}_n}  \frac{|w_-(x,t) - w_-(y,t)|^2}{|x-y|^{N+2s}}\,\dx\dy\dt\\
	&\quad\le
	 \boldsymbol\gm 2^{(N+4)n}\frac{k^2}{\rho^{2s}}|A_n|.
\end{aligned}
\end{equation*}
Departing from the above energy estimate, one can run DeGiorgi's iteration exactly as in \cite[Lemma~3.2]{Liao-cvpd-24} and \cite[Lemma~3.2]{Liao-mod-24}. As a result, we find $\dl$ depending only on the data $\{s, N, C_o, C_1\}$, such that
\begin{equation}\label{Eq:DG:pt-est:2}
u \ge\tfrac12 k \quad
	\mbox{a.e.~in $ K_{\frac12\rho} \times\big(0,\dl\rho^{2s}\big]$.}
\end{equation}
Once $\dl$ is chosen, this also fixes $\nu_o$ from \eqref{Eq:nu:0} in terms of the data $\{s, N, C_o, C_1\}$.

Therefore, no matter which case of \eqref{Eq:DG-alt:1} and \eqref{Eq:DG-alt:2} occurs, one always obtains the desired pointwise estimate from \eqref{Eq:DG:pt-est:1} and \eqref{Eq:DG:pt-est:2}. 
\end{proof}

\begin{remark}\upshape
In Lemma~\ref{Lm:DG:initial:1} a time gap inevitably appears before one could claim pointwise estimate. Technically, this is due to the first case considered in the proof. Whereas the second case is similar to \cite[Lemma~3.2]{Liao-mod-24}.
\end{remark}

\subsection{A measure propagation lemma}
The following measure propagation lemma can be found in \cite[Lemma~3.3]{Liao-mod-24}.

\begin{lemma}\label{Lm:3:1}
 Let $u$ be a   local weak super-solution to \eqref{Eq:1:1}  in $E_T$, such that $u\ge0$.
Introduce parameters $k>0$ and $\al\in (0,1]$. There exist constants $\dl,\,\varep\in(0,1)$ depending only on the data $\{s, N, C_o, C_1\}$ and $\al$, such that if
	\begin{equation*}
	\big|\big\{
		u(\cdot, t_o)\ge k
		\big\}\cap K_{\varrho}(x_o)\big|
		\ge\al |K_{\varrho}|,
	\end{equation*}
	then 
	\begin{equation*}
	\big|\big\{
	u(\cdot, t)\ge \varep k\big\} \cap K_{\varrho}(x_o)\big|
	\ge\frac{\al}2 |K_\varrho|
	\quad\mbox{ for all $t\in\big(t_o,t_o+\dl \varrho^{2s}\big]$,}
\end{equation*}
provided that $K_{\rho}(x_o)\times(t_o,t_o+\dl\rho^{2s}]\subset E_T$. Moreover, we have $\varep=\tfrac18\al$ and $\dl=\bar\dl\al^{3+N}$ for some $\bar\dl\in(0,1)$ depending only on the data $\{s, N, C_o, C_1\}$.
\end{lemma}

\subsection{Measure shrinking lemmata}
The following three lemmata explore the measure shrinking properties of super-solutions. The first one is taken from \cite[Lemma~3.4]{Liao-mod-24}.
\begin{lemma}\label{Lm:3:2-}
 Let $u$ be a  local weak super-solution to \eqref{Eq:1:1}  in $E_T$, such that $u\ge0$.
Introduce parameters $\theta,\, k>0$. Suppose that
	\begin{equation*}
	\big|\big\{
		u(\cdot, t)\ge k
		\big\}\cap K_{\varrho}(x_o)\big|
		\ge\al \big|K_{\varrho}\big|\quad\mbox{ for all $t\in\big(t_o-\theta\varrho^{2s}, t_o\big]$.}
	\end{equation*}
There exists
 $\boldsymbol \gm>1$ depending only on the data $\{s, N, C_o, C_1\}$,  such that 
\begin{equation*}
	\big|\big\{
	u\le \sig k \big\}\cap Q_{\rho}(\theta)\big|
	\le \boldsymbol\gm \frac{\sig }{\theta\al} |Q_{\rho}(\theta)|,
\end{equation*}
provided that $Q_{2\rho}(\theta)\subset E_T$.
\end{lemma}

The previous lemma requires measure theoretical information on each time level and concludes a measure shrinking estimate over a cylinder.
The next lemma weakens the measure theoretical condition, which now is encoded in a local integral, and rebalances the conclusion. As such, it represents the second new element of this section.
\begin{lemma}\label{Lm:3:2}
 Let $u$ be a   local weak super-solution to \eqref{Eq:1:1}  in $E_T$, such that $u\ge0$. Let $k>0$.
There exists a constant
 $\boldsymbol \gm>1$ depending only on the data $\{s, N, C_o, C_1\}$, such that
\begin{equation*}
	\essinf_{t\in[t_o-\rho^{2s}, t_o]}\big|\big\{
	u(\cdot, t)\le  k \big\}\cap  K_{\rho}(x_o)\big|
	\le \frac{ \boldsymbol\gm k }{[u]_{Q_\rho} } |K_{\rho}|,
\end{equation*}
where $[\cdot]_{Q_\rho}$ is the integral average  on $(x_o,t_o)+Q_{\rho}$,
provided that $(x_o,t_o)+Q_{2\rho}\subset E_T$.
\end{lemma}
\begin{proof}
Assume $(x_o,t_o)=(0,0)$. 
Let us  employ the energy estimate \eqref{Eq:cor:2} of Proposition~\ref{Cor:1} in $K_{\rho}\times (-\rho^{2s}, 0]\subset K_{2\rho}\times (-\rho^{2s}, 0]$ 
with the truncation
	$w_- =(u-k)_-$.
Then, we obtain  that
\begin{align*}
	\iint_{Q_\rho}& w_{-}(y,t) \,\dy\dt \bigg(\int_{ \rn }  \frac{w_+(x,t)}{|x-y|^{N+2s}}\,\dx\bigg)\\
	& \le \int_{K_{2\varrho} } w_-^2(x,- \rho^{2s}) \,\dx +
	\frac{\boldsymbol\gm}{\rho^{2s}} \int_{-\rho^{2s}}^{0}\int_{K_{2\rho}} w_{-}^2(x,t) \,\dx \dt\\
	&\quad+\boldsymbol\gm\int_{-\rho^{2s}}^0\int_{K_{2\rho}}  w_{-}(x,t)\,\dx\dt \bigg(  \int_{\rn\setminus K_{\frac32\rho}}\frac{ w_{-}(y,t)}{|y|^{N+2s}}\,\dy\bigg).
\end{align*}
As for the second term on the right-hand side, we use $u\ge0$ and estimate
\begin{align*}
\int_{- \rho^{2s}}^0\int_{K_{2\rho}} & w_{-}(x,t)\,\dx\dt \bigg(  \int_{\rn\setminus K_{\frac32\rho}}\frac{ w_{-}(y,t)}{|y|^{N+2s}}\,\dy\bigg)\\
&\le \boldsymbol\gm k |K_{2\rho} | \bigg(\int_{- \rho^{2s}}^0 \int_{\rn\setminus K_{\frac32\rho}}\frac{ w_{-}(y,t)}{|y|^{N+2s}}\,\dy\bigg)\\
&\le \boldsymbol\gm k^2 |K_{2\rho}|=\boldsymbol\gm 2^N k^2 |K_\rho|.
\end{align*}
The third term is standard:
\[
\int_{K_{2\varrho} } w_-^2(x,- \rho^{2s}) \,\dx\le k^2|K_{2\rho}|=  2^N k^2 |K_\rho|.
\]
Combining the above estimate we see that the energy estimate becomes
\begin{equation}\label{Eq:shrink:0}
	\iint_{Q_\rho} w_{-}(y,t) \,\dy\dt \bigg(\int_{ \rn }  \frac{w_+(x,t)}{|x-y|^{N+2s}}\,\dx\bigg)
	 \le
\boldsymbol\gm k^2 |K_\rho|
\end{equation}
for some $\boldsymbol\gm$ depending only on the data $\{s, N, C_o, C_1\}$.

To estimate the left-hand side of \eqref{Eq:shrink:0}, we evaluate the global integral in the parenthesis only over $K_\rho$ and use that $|x-y|\le 2\rho$ if $x,\,y\in K_\rho$.
As a result, we have that
\begin{align*}
\iint_{Q_\rho }&w_{-}(y,t) \,\dy\dt \bigg(\int_{ K_{\rho}}  \frac{w_+(x,t) }{|x-y|^{N+2s}}\,\dx\bigg)\\
&\ge  \iint_{Q_\rho }w_{-}(y,t) \,\dy\dt \bigg(\int_{ K_{\rho}}  \frac{w_+(x,t) }{(2\rho)^{N+2s}}\,\dx\bigg)\\
&\ge  \essinf_{t\in[-\rho^{2s}, 0]}\int_{K_\rho}w_{-}(y,t) \,\dy \bigg(\iint_{ Q_{\rho}}  \frac{w_+(x,t) }{(2\rho)^{N+2s}}\,\dx\dt \bigg)\\
&\ge   \essinf_{t\in[-\rho^{2s}, 0]}\int_{K_\rho}w_{-}(y,t) \,\dy \bigg(\iint_{ Q_{\rho}}\frac{u(x,t) - k}{(2\rho)^{N+2s}}\,\dx\dt\bigg)\\
&\ge  \essinf_{t\in[-\rho^{2s}, 0]} \int_{K_\rho}w_{-}(y,t) \,\dy \Big( c \biint_{Q_\rho} u\,\dx\dt - c k \Big)\\
&\ge  c \essinf_{t\in[-\rho^{2s}, 0]}\int_{K_\rho}w_{-}(y,t) \,\dy \cdot [u]_{Q_\rho} - c k^2 |K_{\rho}|,
\end{align*}
where $c=c(N)$ is a positive constant, whereas the last integral can be estimated from below by integrating over a smaller set, i.e.,
\[
\essinf_{t\in[-\rho^{2s}, 0]}\int_{K_\rho}w_{-}(y,t) \,\dy \ge \tfrac12 k \essinf_{t\in[-\rho^{2s}, 0]}\big|\big\{u(\cdot, t)\le\tfrac12 k\big\}\cap K_\rho\big|.
\]
Combining these estimates in the energy estimate \eqref{Eq:shrink:0}, we obtain that
\begin{align*}
\essinf_{t\in[-\rho^{2s}, 0]}\big|\big\{u(\cdot, t)\le\tfrac12k\big\}\cap K_\rho\big|\le  \frac{ \boldsymbol\gm k}{[u]_{Q_\rho} } |K_\rho|.
\end{align*}
Redefine $\tfrac12 k$ as $k$ and $2 \boldsymbol\gm$ as $ \boldsymbol\gm$ to conclude.
\end{proof}

Now, we introduce the third  new element of the section. The structure of the statement parallels the last lemma, yet examining the effect of the long-range behavior of $u_+$. 
\begin{lemma}\label{Lm:3:2+}
 Let $u$ be a  local weak super-solution to \eqref{Eq:1:1}  in $E_T$, such that $u\ge0$.
 Let $k>0$ be a parameter. 
There exists a constant
 $\boldsymbol \gm>1$ depending only on the data $\{s, N, C_o, C_1\}$,  such that
\begin{equation*}
\essinf_{t\in[t_o-\rho^{2s}, t_o]}\big|\big\{u(\cdot,t)\le  k \big\}\cap K_\rho(x_o)\big|\le  \frac{ \boldsymbol\gm k}{{\rm Tail}[u_+; (x_o,t_o)+Q_{\rho}]} |K_\rho|,
\end{equation*}
provided that $(x_o,t_o)+Q_{2\rho}\subset E_T$.
\end{lemma}
\begin{proof}
Assume $(x_o,t_o)=(0,0)$. The argument departs from \eqref{Eq:shrink:0}.
To estimate the left-hand side, we integrate instead over $\rn\setminus K_{\rho}$ and use the fact that when $|y|\le \rho$ and $|x|\ge \rho$, one has $|x-y|\le 2 |x|$. As a result, we have that
\begin{align*}
\iint_{Q_\rho }&w_{-}(y,t) \,\dy\dt \bigg(\int_{\rn\setminus K_{\rho}}  \frac{w_+(x,t) }{|x-y|^{N+2s}}\,\dx\bigg)\\
&\ge 2^{-(N+2s)}\iint_{Q_\rho }w_{-}(y,t) \,\dy\dt \bigg(\int_{\rn\setminus K_{\rho}}  \frac{w_+(x,t) }{|x|^{N+2s}}\,\dx\bigg)\\
&\ge 2^{-(N+2s)}\essinf_{t\in[-\rho^{2s}, 0]}\int_{K_\rho}w_{-}(y,t) \,\dy \bigg(\int_{- \rho^{2s}}^0\int_{\rn\setminus K_{\rho}}  \frac{w_+(x,t) }{|x|^{N+2s}}\,\dx\dt \bigg)\\
&\ge 2^{-(N+2s)} \essinf_{t\in[-\rho^{2s}, 0]}\int_{K_\rho}w_{-}(y,t) \,\dy \bigg(\int_{- \rho^{2s}}^0\int_{\rn\setminus K_{\rho}}\frac{u_+(x,t) - k}{|x|^{N+2s}}\,\dx\dt\bigg)\\
&\ge 2^{-(N+2s)}\essinf_{t\in[-\rho^{2s}, 0]} \int_{K_\rho}w_{-}(y,t) \,\dy \Big({\rm Tail}[u_+, Q_{\rho}] - \boldsymbol\gm k \Big)\\
&\ge 2^{-(N+2s)} \essinf_{t\in[-\rho^{2s}, 0]}\int_{K_\rho}w_{-}(y,t) \,\dy \cdot {\rm Tail}[u_+, Q_{\rho}] - \boldsymbol\gm k^2 |K_{\rho}|.
\end{align*}
The proof is then concluded just like in Lemma~\ref{Lm:3:2}.
\end{proof}


\section{Nonlocal parabolic weak Harnack estimates}\label{S:Harnack}
This section is devoted to the proof of nonlocal parabolic weak Harnack estimates.
The structure is as follows. First in \S~\ref{S:exp-pos} we use Lemma~\ref{Lm:DG:1}, Lemma~\ref{Lm:3:1} and Lemma~\ref{Lm:3:2-} to show an expansion of positivity (Proposition~\ref{Prop:expansion}), from which Theorem~\ref{Prop:WHI:1} follows in \S~\ref{S:WHI-proof:1}. Then, we use Lemma~\ref{Lm:DG:initial:1} and Lemma~\ref{Lm:3:2+} (Lemma~\ref{Lm:3:2}, respectively), with the aid of Proposition~\ref{Prop:expansion} to show Theorem~\ref{Prop:WHI:2} (Theorem~\ref{Prop:WHI:2-}, respectively) in \S~\ref{S:WHI-proof:2}.
\subsection{Expansion of positivity}\label{S:exp-pos}

\begin{proposition}\label{Prop:expansion}
Let $u$ be a  local weak super-solution to \eqref{Eq:1:1}   in $E_T$, such that $u\ge0$.
Suppose for some constants  $\al\in(0,1]$ and $k>0$, there holds
	\begin{equation*}
		\big|\big\{u(\cdot, t_o)\ge k \big\}\cap K_{\varrho}(x_o) \big|
		\ge
		\al |K_\varrho |.
	\end{equation*}
There exists a constant   
 $\eta\in(0,1)$ depending on the data $\{s, N, C_o, C_1\}$ and $\al$, such that 
\begin{equation*}
	u\ge\eta k
	\quad
	\mbox{a.e.~in $ K_{2\varrho}(x_o) \times\big( t_o+\tfrac12 (8\varrho)^{2s},
	t_o+2(8\varrho)^{2s}\big],$}
\end{equation*}
provided  
\[
K_{4\rho}(x_o)\times\big(t_o, t_o+4(8\varrho)^{2s}\big]\subset E_T.
\]
Moreover, $\eta=\bar{\eta}\al^q$ for some $\bar\eta\in(0,1)$ and $q>1$ depending on the data $\{s, N, C_o, C_1\}$.
\end{proposition}
\begin{proof}
Assume $(x_o,t_o)=(0,0)$ for simplicity.
 Rewriting the measure theoretical information at the initial time $t_o=0$ in the larger ball $K_{4\rho}$ and replacing $\al$ by $4^{-N}\al$, we 
 apply Lemma~\ref{Lm:3:1} to obtain $\varep=\tfrac18\al$ and $\dl=\bar{\dl}\al^{3+N}$ for $\bar\dl\in(0,1)$ depending only on the data $\{s,  N, C_o, C_1\}$, such that 
	\begin{equation}\label{Eq:exp-pos:1}
	\big|\big\{
	 u(\cdot, t) \ge \varep k\big\} \cap K_{4\varrho} \big|
	\ge\frac{\al}2 4^{-N} |K_{4\varrho}|
	\quad\mbox{ for all $t\in\big(0, \dl (4\varrho)^{2s}\big]$.}
\end{equation}

Notice that 
\[
(0,\bar{t})+Q_{4\rho}(\tfrac12\dl )\subset K_{4\rho}\times\big(0, \dl (4\varrho)^{2s}\big]
\]
whenever
\begin{equation}\label{Eq:t-bar}
\bar{t}\in \big(\tfrac12\dl (4\varrho)^{2s}, \dl (4\varrho)^{2s}\big].
\end{equation}
Hence, the measure theoretical information \eqref{Eq:exp-pos:1}  allows us to apply Lemma~\ref{Lm:3:2-} in the cylinders $(0,\bar{t})+Q_{4\rho}(\tfrac12\dl )$ as $\bar{t}$ ranges over the time interval in \eqref{Eq:t-bar},
and with $k$ and $\al$ replaced by $\varep k$ and $\tfrac12 4^{-N}\al$. 

Letting $\nu=\bar\nu\dl^{q}$ be determined in Lemma~\ref{Lm:DG:1} for some $\bar\nu$ and $q$  in terms of the data, we further choose $\sig$ according to Lemma~\ref{Lm:3:2-} to satisfy
\[
\boldsymbol\gm \frac{\sig }{\dl\al} \le \nu, \quad\text{i.e.}\quad \sig=\frac{\nu\dl\al}{\boldsymbol\gm}=\frac{\bar\nu\bar\dl^{q+1}\al^{(3+N)(q+1)}}{\boldsymbol\gm},
\]
where we used $\dl=\bar{\dl}\al^{3+N}$ to calculate.
The above choice is viable as $\boldsymbol\gm$ of Lemma~\ref{Lm:3:2-} is independent of $\sig$.
Letting $\widetilde{\boldsymbol\gm}$ be chosen in Lemma~\ref{Lm:DG:1} 
such a choice of $\sig$ permits us to first apply Lemma~\ref{Lm:3:2-} (with $\varep k$ and $\tfrac12 4^{-N}\al$) and then Lemma~\ref{Lm:DG:1} (with $\sig\varep k$) in the cylinder $(0,\bar{t})+Q_{4\rho}(\tfrac12\dl)$ with an arbitrary $\bar{t}$ as in \eqref{Eq:t-bar}. Therefore, by arbitrariness of $\bar{t}$ we conclude that
\[
u\ge\tfrac14\sig\varep k \quad\text{a.e. in}\>K_{2\varrho}  \times\big(\tfrac12 \dl (4\rho)^{2s}, \dl  (4\varrho)^{2s}\big].
\]

Once pointwise estimate is reached, we can repeat the above procedure with $\al=1$.
In fact, let us denote $\eta=\tfrac14\sig\varep$, and use the above pointwise estimate to run the procedure and obtain
\[
u\ge\bar{\eta} \eta k \quad\text{a.e. in}\>K_{2\varrho}  \times\big(\tfrac12 \dl (4\rho)^{2s}+\tfrac12 \bar{\dl} (8\rho)^{2s}, \dl (4\rho)^{2s}+ \bar{\dl}  (8\varrho)^{2s}\big]
\]
for some $\bar{\eta},\,\bar{\dl}\in(0,\frac12)$ depending only on the data and independent of $\al$. 
Repeating this procedure we obtain for $n\in\nn$ that 
\begin{equation}\label{Eq:pt-est-n}
u\ge\bar{\eta}^n \eta k \quad\text{a.e. in}\>K_{2\varrho}  \times\big(\tfrac12 \dl (4\rho)^{2s}+\tfrac12 \bar{\dl} (8\rho)^{2s}, \dl (4\rho)^{2s}+ n\bar{\dl}  (8\varrho)^{2s}\big]
\end{equation}
for the same $\bar{\eta},\,\bar{\dl}\in(0,\frac12)$.

Finally, we only need to select $n=\lceil 2/\bar{\dl}\rceil$. In this way, the pointwise estimate \eqref{Eq:pt-est-n} yields that
\[
u\ge\bar\eta^{\lceil 2/\bar{\dl}\rceil} \eta k \quad\text{a.e. in}\>K_{2\varrho}  \times\big( \bar{\dl} (8\rho)^{2s},  2(8\varrho)^{2s}\big].
\]
Recall the dependence of $\dl$ and $\eta$ on $\al$, and redefine relevant constants to conclude.
\end{proof}
\begin{remark}\upshape
The key information provided by Proposition~\ref{Prop:expansion} is that $\eta$ has a power-like dependence on $\al$.
In the case of local operators, such a result often requires more work to reach. See~\S~\ref{N-S-p} for more discussions.
\end{remark}
\subsection{Proof of Theorem~\ref{Prop:WHI:1}}\label{S:WHI-proof:1}
Assume $(x_o,t_o)=(0,0)$ and define
\[
{\bf I}:=  \essinf_{K_{2\rho}\times(\frac12 (8\varrho)^{2s}, 2(8\varrho)^{2s}]}u.
\]
For $\varep>0$ to be fixed, let us first estimate the $L^\varep$-norm of $u(\cdot,0)$ by its measure distribution:
\begin{equation}\label{Eq:7:4}
\begin{aligned}
\int_{K_{\rho}}u^{\varep}(\cdot,0)\,\dx&= \varep\int_0^\infty |[u(\cdot,0)>k]\cap K_\rho | k^{\varep-1}\,\d k\\
&\le \varep\int_{\bf I}^\infty |[u(\cdot,0)>k]\cap K_\rho | k^{\varep-1}\,\d k+{\bf I}^{\varep}|K_{\rho}|.
\end{aligned}
\end{equation}
By Proposition~\ref{Prop:expansion}, there exist $q>1$ and $\eta\in(0,1)$
depending only on the data $\{s, N, C_o, C_1\}$, such that
\[
{\bf I}\ge\eta k\bigg(\frac{|[u(\cdot,0)>k]\cap K_{\rho}|}{|K_\rho|}\bigg)^q.
\]
Using this we estimate the first term on the right-hand side of \eqref{Eq:7:4} by
\[
\varep\int_{\bf I}^\infty |[u(\cdot,0)>k]\cap K_\rho|k^{\varep-1}\,\d k\le\frac{\varep {\bf I}^{\frac1q}}{\eta^{\frac1q}}|K_{\rho}|\int_{\bf I}^\infty k^{\varep-\frac1q-1}\,\d k.
\]
Now, we stipulate to take $\varep<1/q$, such that the improper integral on the right-hand side
converges. In such a way, the right-hand side of \eqref{Eq:7:4} is 
bounded by 
$$
\bigg[1+\frac{\varep}{(\varep-\frac1q)\eta^{\frac1q}}\bigg] {\bf I}^\varep |K_\rho|.
$$ 
Substituting it back in \eqref{Eq:7:4},
 we obtain the desired conclusion.
\subsection{Proof of Theorems~\ref{Prop:WHI:2-}~\&~\ref{Prop:WHI:2}}\label{S:WHI-proof:2}
Since the arguments are similar, we only show the proof of Theorem~\ref{Prop:WHI:2} in detail.

Let $\nu_o$ be fixed in Lemma~\ref{Lm:DG:initial:1} in terms of the data $\{s, N, C_o, C_1\}$. We apply Lemma~\ref{Lm:3:2+} and choose $k$ to satisfy
\begin{equation}\label{Eq:exp2:1}
\essinf_{t\in[-(2\rho)^{2s}, 0]}\big|\big\{u(\cdot,t)\le  k \big\}\cap K_{2\rho}\big|\le \frac{ \boldsymbol\gm k}{{\rm Tail}[u_+, Q_{2\rho}]} |K_{2\rho}|
\le \nu_o|K_{2\rho}|.
\end{equation}
That is, we choose
\begin{equation}\label{Eq:choice-k}
k=\frac{\nu_o}{\boldsymbol\gm}{\rm Tail}[u_+, Q_{2\rho}].
\end{equation}
Let $\bar{t}\in[-(2\rho)^{2s}, 0]$ be an instant when the infimum of \eqref{Eq:exp2:1} is attained. Then, according to Lemma~\ref{Lm:DG:initial:1} we have for some $\dl\in(0,\frac14)$
\begin{equation}\label{Eq:exp2:2}
u\ge\tfrac1{16}k \quad\mbox{a.e.~in $K_{\frac12\rho}\times\big(\bar{t}+\tfrac34\dl \rho^{2s}, \bar{t}+\dl \rho^{2s}\big]$,}
\end{equation}
provided that $K_{2\rho}\times[-(2\rho)^{2s}, \dl\rho^{2s}]\subset \mathcal{Q}$, 
where $\dl$ and $\widetilde{\boldsymbol\gm}$ depend only on the data $\{s, N, C_o, C_1\}$.

Next, we use \eqref{Eq:exp2:2} and apply Proposition~\ref{Prop:expansion} with $\al=1$. Consequently, we obtain
some $\eta$ depending only on the data $\{s, N, C_o, C_1\}$, such that
\begin{equation}\label{Eq:exp2:3}
	u\ge\eta k
	\quad
	\mbox{a.e.~in $ K_{\varrho}  \times\big( \bar{t} + \tfrac34\dl \rho^{2s} +\tfrac12 (4\varrho)^{2s},
	\bar{t}+\dl \rho^{2s}+2(4\varrho)^{2s}\big],$}
\end{equation}
provided that $K_{2\rho}\times(-(2\rho)^{2s}, 6(4\varrho)^{2s}]\subset \mathcal{Q}$.
No matter where $\bar{t}$ is in $[-(2\rho)^{2s}, 0]$, we can always conclude from \eqref{Eq:exp2:3} that
\begin{equation*}
	u\ge\eta k
	\quad
	\mbox{a.e.~in $ K_{\varrho}  \times\big( \tfrac34 (4\varrho)^{2s},
	(4\varrho)^{2s}\big],$}
\end{equation*}
provided that  $K_{2\rho}\times(-(2\varrho)^{2s}, 6(4\varrho)^{2s}]\subset \mathcal{Q}$.
The proof is concluded by redefining $\eta\nu_o/(\widetilde{\boldsymbol\gm}\boldsymbol\gm)$ as $\eta$ and recall $k$ from \eqref{Eq:choice-k}.

As for the dimensionless version Theorem~\ref{Prop:WHI:3}, one only needs to observe that all lemmata in \S~\ref{S:prelim} can be produced with $\rho=1$ under the assumption that $u$ solves \eqref{Eq:1:1} in $K_4\times(0,4]$ and 
\[
\|b\|_{L^{\infty}(0,4;L^{\frac{N}s}(K_2))}\le C_3.
\]
The generic constant $\boldsymbol\gm$ now depends on the data $\{s, N, C_o, C_1\}$ and $C_3$.
Consequently, the expansion of positivity also admits a dimensionless version which eventually along the same lines of reasoning yields Theorem~\ref{Prop:WHI:3}.

\section{$L^{\infty}_{\loc}$--$L^{\varep}_{\loc}$ estimates for sub-solutions}\label{S:bd}

\subsection{Proof of Theorem~\ref{Prop:A:1}}
The proof hinges upon the energy estimate in Proposition~\ref{Prop:2:1}.
First suppose $Q(R,\theta)\subset E_T$, and let $\rho\in(0,R)$.
Let us assume $(x_o,t_o)=(0,0)$ for simplicity. 
For $n\in\nn_0$, $\sig\in(0,1)$, and $k>0$, introduce 
\begin{align*}
	\left\{
	\begin{array}{c}
 k_n= k - 2^{-n} k, 
\\[5pt]
 \rho_n=\sig\rho + 2^{-n}(1-\sig)\rho,\quad  \theta_n=\sig \theta + 2^{-n}(1-\sig)\theta ,\\[5pt]
 \widetilde{\rho}_n=\frac12(\rho_n+\rho_{n+1}), \quad \widehat{\rho}_n=\frac34\rho_n+ \frac14\rho_{n+1},\quad \overline{\rho}_n=\frac14\rho_n+ \frac34\rho_{n+1},\\[5pt]
 \widetilde{\theta}_n=\frac12(\theta_n+\theta_{n+1}), \quad \widehat{\theta}_n=\frac34\theta_n+ \frac14\theta_{n+1},\quad \overline{\theta}_n=\frac14\theta_n+ \frac34\theta_{n+1},\\[5pt]
  K_n=K_{\rho_n}, \quad \widetilde{K}_n=K_{\widetilde{\rho}_n},\quad \widehat{K}_n=K_{\widehat\rho_n},\quad \overline{K}_n=K_{\overline{\rho}_n},\\[5pt]
\dsty Q_n=K_n\times(-\theta_n,0),\quad\widetilde{Q}_n=\widetilde{K}_n\times(-\widetilde{\theta}_n,0),\\[5pt]
\widehat{Q}_n=\widehat{K}_n\times(-\widehat{\theta}_n,0),\quad\overline{Q}_n=\overline{K}_n\times(-\overline{\theta}_n,0).
\end{array}
	\right.
\end{align*}
Observe that $Q_o=Q(\rho,\theta)$, $Q_{\infty}=Q(\sig\rho,\sig\theta)$, and $Q_{n+1}\subset\overline{Q}_n\subset\widetilde{Q}_n\subset\widehat{Q}_n\subset Q_n$.
Introduce the cutoff function $\z$ in $Q_n$, vanishing outside $\widehat{Q}_{n}$ and
equal to the identity in $\widetilde{Q}_{n}$, such that
\begin{equation*}
|D\z|\le \frac{2^{n+4}}{(1-\sig)\rho}\quad\text{ and }\quad |\pl_t \z|\le \frac{2^{n+4}}{(1-\sig)\theta}.
\end{equation*}
The energy estimate of Proposition~\ref{Prop:2:1} is used in $Q_n$ with $\z$ and
with 
$$
w_+(x,t):=\big(u(x,t)-\ell(t)-k_{n+1}\big)_+, \quad \ell(t):=\widetilde{\boldsymbol\gm}\int^t_{-\theta}\int_{\rn\setminus K_R}\frac{ u_{+}(y,\tau)}{|y|^{N+2s}}\,\dy\d\tau,
$$
where $\widetilde{\boldsymbol\gm}>0$ is to be determined.
As a result, we have
\begin{align}\nonumber
	\essup_{t\in[-\widetilde{\theta}_n,0]}&\int_{\widetilde{K}_n} w^2_{+}(x,t)\,\dx+ \int_{-\widetilde{\theta}_n}^{0}\int_{\widetilde{K}_n}\int_{\widetilde{K}_n}  \frac{|w_+(x,t) - w_+(y,t)|^2}{|x-y|^{N+2s}}\,\dx\dy\dt\\\nonumber
	&\quad\le
	\boldsymbol\gm\int_{-\theta_n}^{0}\int_{K_n}\int_{K_n}\max\big\{w^2_{+}(x,t), w^2_{+}(y,t)\big\} \frac{|\z(x,t) - \z(y,t)|^2}{|x-y|^{N+2s}}\,\dx\dy\dt\\\nonumber
	&\qquad+\boldsymbol\gm\iint_{Q_n} \z^2 w_{+}(x,t)\,\dx\dt \bigg(\essup_{\substack{x\in \widehat{K}_n}} \int_{\rn\setminus K_n}\frac{ w_{+}(y,t)}{|x-y|^{N+2s}}\,\dy\bigg)\\\label{Eq:energy-DG}
	&\qquad-2\iint_{Q_n}  \ell^{\prime}(t) \z^2 w_{+}(x,t)\,\dx\dt + \iint_{Q_n} |\pl_t\z^2|w_{+}^2(x,t)\,\dx\dt.
\end{align}

The four terms on the right-hand side of \eqref{Eq:energy-DG} are treated as follows.
The first term is standard. Indeed, we estimate
\begin{align*}
\int_{-\theta_n}^{0}&\int_{K_n}\int_{K_n}\max\big\{w^2_{+}(x,t), w^2_{+}(y,t)\big\} \frac{|\z(x,t) - \z(y,t)|^2}{|x-y|^{N+2s}}\,\dx\dy\dt\\
&\le  \frac{2^{2n+9}}{(1-\sig)^2\rho^{2}}\int_{-\theta_n}^{0}\int_{K_n}\int_{K_n} \frac{w^2_{+}(x,t)}{|x-y|^{N+2(s-1)}}\,\dx\dy\dt\\
&\le   \frac{\boldsymbol \gm 2^{2n}}{(1-\sig)^2\rho^{2s}} \iint_{Q_n} w^2_{+}(x,t)\,\dx\dt.
\end{align*}
The last term is also standard, namely,
\[
\iint_{Q_n} |\pl_t\z^2|w_{+}^2(x,t)\,\dx\dt \le\frac{2^{n+4}}{(1-\sig)\theta } \iint_{Q_n} w_{+}^2(x,t)\,\dx\dt.
\]
 
The second term and the third, negative term need to be packed.  To this end, observe that when $|y|\ge \rho_n$ and $|x|\le \widehat\rho_n$, there holds
\[
\frac{|y-x|}{|y|}\ge1-\frac{\widehat\rho_n}{ \rho_n}=\frac14\Big(\frac{\rho_n-\rho_{n+1}}{\rho_n}\Big)\ge\frac{1-\sig}{2^{n+3}};
\]
when $|y|\ge R$ and $|x|\le \rho$, there holds
\[
\frac{|y-x|}{|y|}\ge1-\frac{|x|}{|y|}\ge 1-\frac{\rho}{ R}. 
\]
Consequently, 
we estimate the second term as
\begin{align}\nonumber
\iint_{Q_n}&  \z^2 w_{+}(x,t)\,\dx\dt \bigg[\essup_{\substack{x\in\widehat{K}_n}} \int_{\rn\setminus K_n}\frac{ w_{+}(y,t)}{|x-y|^{N+2s}}\,\dy\bigg]\\ \nonumber
&=\iint_{Q_n}  \z^2 w_{+}(x,t)\,\dx\dt \essup_{\substack{x\in\widehat{K}_n}} \bigg[\int_{K_R\setminus K_n}\frac{ w_{+}(y,t)}{|x-y|^{N+2s}}\,\dy+\int_{\rn\setminus K_R}\frac{ w_{+}(y,t)}{|x-y|^{N+2s}}\,\dy\bigg]\\ \nonumber
&\le   \frac{\boldsymbol\gm 2^{(N+2s)n}}{(1-\sig)^{N+2s}} \iint_{Q_n}  \z^2 w_{+}(x,t)\,\dx\dt \bigg[ \int_{K_R\setminus K_n}\frac{ w_{+}(y,t)}{|y|^{N+2s}}\,\dy\bigg]\\ \nonumber
&\qquad+  \Big(\frac{R}{R-\rho}\Big)^{N+2s}   \iint_{Q_n}  \z^2 w_{+}(x,t)\,\dx\dt \bigg[ \int_{\rn\setminus K_R}\frac{ w_{+}(y,t)}{|y|^{N+2s}}\,\dy\bigg]\\ \nonumber
&\le  \frac{\boldsymbol\gm 2^{(N+2s)n}}{[\sig (1-\sig)]^{N+2s}\rho^{2s}}\Big(\frac{R}{\rho}\Big)^{N}\essup_{t\in[-\theta,0]} \bint_{K_R} u_{+}(y,t)\,\dy \iint_{Q_n}  \z^2 w_{+}(x,t)\,\dx\dt \\ \nonumber
&\qquad+\Big(\frac{R}{R-\rho}\Big)^{N+2s}   \iint_{Q_n}  \z^2 w_{+}(x,t)\,\dx\dt  \bigg[ \int_{\rn\setminus K_R}\frac{ u_{+}(y,t)}{|y|^{N+2s}}\,\dy\bigg]. 
\end{align}
The last term in the above estimate will cancel with  the third, negative term on the right-hand side of the energy estimate \eqref{Eq:energy-DG},
if we choose 
$$
2\widetilde{\boldsymbol\gm}=\Big(\frac{R}{R-\rho}\Big)^{N+2s}\quad\implies\quad \ell(t)=\frac{1}{2}\Big(\frac{R}{R-\rho}\Big)^{N+2s}\int^t_{-\theta}\int_{\rn\setminus K_R}\frac{ u_{+}(y,\tau)}{|y|^{N+2s}}\,\dy\d\tau.
$$ As a result of this choice, the second and the third terms in \eqref{Eq:energy-DG} together are bounded by
\[
 \frac{\boldsymbol\gm 2^{(N+2s)n}\mathbf{I}}{[\sig(1-\sig)]^{N+2s}\rho^{2s}} \iint_{Q_n}  \z^2 w_{+}(x,t)\,\dx\dt
 \]
where
 \[
 \mathbf{I}:=\Big(\frac{R}{\rho}\Big)^{N}\essup_{t\in[-\theta,0]} \bint_{K_R} u_{+}(y,t)\,\dy.
\]
To further estimate the last integral, we use H\"older's inequality and the fact that
\begin{equation}\label{Eq:obs}
|A_n|:=\big|\big\{u(x,t) -\ell(t)\ge k_{n+1}\big\}\cap Q_n \big| \le\frac{2^{2(n+1)}}{k^2} \iint_{Q_n}   \widetilde{w}^2_{+}(x,t)\,\dx\dt,
\end{equation}
where we denoted
\[
\widetilde{w}_+(x,t):=\big(u(x,t)-\ell(t)-k_{n}\big)_+.
\]
As a result, noticing also that $w_+\le \widetilde{w}_+$, we have
\begin{align*}
\iint_{Q_n}  \z^2 w_{+}(x,t)\,\dx\dt &\le \Big[\iint_{Q_n}  w^2_{+}(x,t)\,\dx\dt \Big]^{\frac12} \big|\big\{u(x,t) -\ell(t)\ge k_{n+1}\big\}\cap Q_n \big|^{\frac12}\\
&\le \frac{2^{n+1}}{k} \iint_{Q_n}   \widetilde{w}^2_{+}(x,t)\,\dx\dt.
\end{align*}
Collecting these estimates on the right-hand side of \eqref{Eq:energy-DG} and using $w_+\le \widetilde{w}_+$, we arrive at
\begin{align*}
\essup_{t\in[-\widetilde\theta_n,0]}&\int_{\widetilde{K}_n} w_+^2\,\dx +
\int_{-\widetilde\theta_n}^{0}\int_{\widetilde{K}_n}\int_{\widetilde{K}_n}  \frac{|w_+(x,t) - w_+(y,t)|^2}{|x-y|^{N+2s}}\,\dx\dy\dt\\
&\le \frac{\boldsymbol\gm 2^{(N+2)n}}{(1-\sig)^{N+2}}\Big(\frac1{\theta} +\frac1{\rho^{2s}} +\frac1{\sig^{N+2s}\rho^{2s}}\frac{\mathbf{I}}{k}\Big)  \iint_{Q_n}   \widetilde{w}^2_{+}\,\dx\dt\\
&\le  \frac{\boldsymbol\gm 2^{(N+2)n}}{(1-\sig)^{N+2}\rho^{2s}}\Big(\frac{\rho^{2s}}{\theta}+\frac1{\dl\sig^{N+2s}}\Big) \iint_{Q_n}   \widetilde{w}^2_{+}\,\dx\dt.
\end{align*}
To obtain the last line, we enforced that
\begin{equation}\label{Eq:k-M}
k\ge \dl \mathbf{I}.
\end{equation}

Now set $0\le\phi\le1$ to be a cutoff function in $\widetilde{Q}_n$, which vanishes outside $\overline{Q}_n$,
 equals the identity in $Q_{n+1}$ and satisfies $|D\phi|\le 2^{n+4}/\rho$.
 An application of the H\"older inequality  and the Sobolev imbedding
(cf.~\cite[Proposition~A.3]{Liao-cvpd-24} with $d=2^{-n-4}$)  gives that 
\begin{align*}
	&\iint_{Q_{n+1}} w_+^2\,\dx\dt
	\le 
	\iint_{\widetilde{Q}_n}(\phi w_+)^2\,\dx\dt\\
	&\quad\le
	\Big[\iint_{\widetilde{Q}_n}(\phi w_+)^{2\kappa }
	\,\dx\dt\Big]^{\frac{1}{\kappa }}|A_n|^{1-\frac{1}{\kappa }}\\
	&\quad\le\boldsymbol\gm
	\Big[\rho^{2s}\int_{-\widetilde\theta_n}^{0}\iint_{\widetilde{K}_n^2}  \frac{|\phi w_+(x,t) - \phi w_+(y,t)|^2}{|x-y|^{N+2s}}\,\dx\dy\dt+\frac1{d^{N+2s}}\iint_{\widetilde{Q}_n}(\phi w_+)^2\,
	\dx\dt\Big]^{\frac{1}{\kappa }}\\
	&\qquad\ 
	\cdot\Big[\essup_{t\in[-\widetilde\theta_n,0]}
	\bint_{\widetilde{K}_n}(\phi w_+)^2\,\dx\Big]^{\frac{(\kappa_*-1)}{\kappa_* \kappa }}
	 |A_n|^{1-\frac{1}{\kappa }}
\end{align*}
where $|A_n|$ is defined in \eqref{Eq:obs} and $\kappa:=2-1/\kappa_*$ with
\begin{equation*}
\kappa_*:=\left\{
\begin{array}{cl}
\frac{N}{N-2s} \quad&\text{if}\quad 2s<N,\\[5pt]
\in (1,\infty) \quad &\text{if} \quad 2s \ge N.
\end{array}\right.
\end{equation*} 
To estimate the fractional norm, we use the triangle inequality
\begin{align*}
|w_+\phi(x,t) &- w_+\phi(y,t)|^2 \\
&\le 2   |w_+(x,t) - w_+(y,t)  |^2 \phi^2(x,t) 
 + 2  w^2_+(y,t)   |\phi(x,t) - \phi(y,t) |^2,
\end{align*}
such that
\begin{align*}
\rho^{2s}&\int_{-\widetilde\theta_n}^{0}\int_{\widetilde{K}_n}\int_{\widetilde{K}_n}  \frac{|w_+\phi(x,t) - w_+\phi(y,t)|^2}{|x-y|^{N+2s}}\,\dx\dy\dt\\
&\le 2 \rho^{2s}\int_{-\widetilde\theta_n}^{0}\int_{\widetilde{K}_n}\int_{\widetilde{K}_n}  \frac{|w_+(x,t) - w_+(y,t)|^2}{|x-y|^{N+2s}}\,\dx\dy\dt\\
&\quad + 2 \rho^{2s}\int_{-\widetilde\theta_n}^{0}\int_{\widetilde{K}_n}\int_{\widetilde{K}_n}  \frac{w^2_+(y,t)|\phi(x,t) - \phi(y,t)|^2}{|x-y|^{N+2s}}\,\dx\dy\dt\\
&\le 2 \rho^{2s}\int_{-\widetilde\theta_n}^{0}\int_{\widetilde{K}_n}\int_{\widetilde{K}_n}  \frac{|w_+(x,t) - w_+(y,t)|^2}{|x-y|^{N+2s}}\,\dx\dy\dt\\
&\quad + \boldsymbol\gm 2^{2n} \iint_{\widetilde{Q}_n}  w^2_+(y,t) \,\dy\dt\\
&\le  \frac{\boldsymbol\gm 2^{(N+2)n}}{(1-\sig)^{N+2}\rho^{2s}}\Big(\frac{\rho^{2s}}{\theta}+\frac1{\dl\sig^{N+2s}}\Big) \iint_{Q_n}   \widetilde{w}^2_{+}(x,t)\,\dx\dt.
\end{align*}
Here, to obtain the last line we applied the energy estimate, and also used that $w_+\le \widetilde{w}_+$ and $\widetilde{Q}_n\subset Q_n$. A straightforward computation shows that the term containing $d$ can be estimated by the same quantity as in the last display, and  the term with $\essup$ is also estimated by this quantity according to the energy estimate. Whereas the term $|A_n|$ is estimated by \eqref{Eq:obs}.
Therefore, we arrive at
\begin{align*}
\iint_{Q_{n+1}} & w_+^2\,\dx\dt\\
&\le 
	\frac{\dsty\boldsymbol\gm\boldsymbol b^n  \Big(\frac{\rho^{2s}}{\theta}+\frac1{\dl\sig^{N+2s}}\Big)^{\frac{2\kappa_*-1}{\kappa_*\kappa }}}{(1-\sig)^{(N+2)\frac{2\kappa_*-1}{\kappa_*\kappa }}} \rho^{-(N+2s)\frac{\kappa_* -1}{\kappa_*\kappa }} k^{2(\frac1\kappa -1)}  
	\Big(\iint_{Q_{n}} \widetilde{w}_+^2\,\dx\dt\Big)^{1+\frac{\kappa_* -1}{\kappa_*\kappa }}
\end{align*}
for some $\boldsymbol b=\boldsymbol b(N)>1$.
At this stage, we denote
\[
\boldsymbol{Y}_n:=\iint_{Q_{n}} (u(x,t)-\ell(t) - k_{n})_+^2\,\dx\dt,
\]
recall $\kappa=2-\frac1{\kappa_*}$ so that $1-\frac1\kappa=\frac{\kappa_*-1}{\kappa_*\kappa}$, and rewrite the previous recursive estimate as
\[
\boldsymbol{Y}_{n+1}\le \frac{\dsty\boldsymbol\gm\boldsymbol b^n\Big(\frac{\rho^{2s}}{\theta}+\frac1{\dl\sig^{N+2s}}\Big)^{\frac{2\kappa_*-1}{\kappa_*\kappa }}}{(1-\sig)^{(N+2)\frac{2\kappa_*-1}{\kappa_*\kappa }}} 
\frac{\boldsymbol{Y}_n^{1+\frac{\kappa_* -1}{\kappa_*\kappa }}}{\dsty (\rho^{N+2s} k^2)^{\frac{\kappa_* -1}{\kappa_*\kappa }}}.
\]
Hence, by the fast geometric convergence, cf. \cite[Chapter I, Lemma~4.1]{DB}, 
there exists
a  constant $\boldsymbol\gm>1$ depending only on the data $\{s, N, C_o, C_1\}$, such that
if we require  
\begin{equation}\label{Eq:k-M:1}
\boldsymbol{Y}_o=\iint_{Q_{o}} (u-\ell)_+^2\,\dx\dt\le  \frac{(1-\sig)^{\frac{(N+2)(2\kappa_*-1)}{\kappa_* -1 }}}{\dsty\boldsymbol\gm\Big(\frac{\rho^{2s}}{\theta}+\frac1{\dl\sig^{N+2s}}\Big)^{\frac{2\kappa_*-1}{\kappa_*-1 }}} \rho^{N+2s} k^2,
\end{equation}
then
$$
\lim_{n\to\infty}\iint_{Q_{n}} (u-\ell -k_n)_+^2\,\dx\dt=0, \quad\text{i.e.}\quad u(x,t) \le \ell(t) + k\quad\text{a.e. in}\> Q_\infty=Q(\sig\rho,\sig\theta).
$$ 
Taking both requirements \eqref{Eq:k-M} and \eqref{Eq:k-M:1}  of $k$ into account, we obtain that
\begin{align}\nonumber
\essup_{Q(\sig\rho,\sig \theta)} u &\le \dl  \mathbf{I} + \frac{1}{2}\Big(\frac{R}{R-\rho}\Big)^{N+2s}\int^0_{-\theta}\int_{\rn\setminus K_R}\frac{ u_{+}(y,\tau)}{|y|^{N+2s}}\,\dy\d\tau\\ \label{Eq:A:4}
&\quad+ \boldsymbol\gm (1-\sig)^{-(N+2)q}
\underbrace{\Big(\frac{\rho^{2s}}{\theta}+\frac1{\dl\sig^{N+2s}}\Big)^{q}\Big(\frac{\theta}{\rho^{2s}} \Big)^{\frac12}  }_{:=\mathcal{A}_\dl}\Big(\biint_{Q(\rho, \theta)} u_+^2\,\dx\dt\Big)^{\frac12}.
\end{align}
For simplicity, we denoted $q:=\frac{2\kappa_*-1}{2(\kappa_*-1) }$ in \eqref{Eq:A:4}.
The right-hand side of \eqref{Eq:A:4} is finite due to the notion of solution, and hence $u_+$ is  locally bounded in $E_T$.

Now, we refine \eqref{Eq:A:4}, which holds for all $Q(R,\theta)\subset E_T$, $\dl,\,\sig\in(0,1)$ and $\rho\in(0,R)$.
The task is twofold: To reduce the integral exponent in the last term of  \eqref{Eq:A:4} to $\varep$; and, to absorb the first term into the left-hand side. In order for that, let $\varep\in(0,2)$. The last term in \eqref{Eq:A:4} is estimated by Young's inequality:
\begin{align*}
&\frac{\boldsymbol\gm \mathcal{A}_\dl}{ (1-\sig)^{(N+2)q-\frac{N+1}{2}}} \Big[\frac1{(1-\sig)^{N+1}}\biint_{Q(\rho, \theta)} u_+^2\,\dx\dt\Big]^{\frac12}\\
&\qquad\le\frac{\boldsymbol\gm \mathcal{A}_\dl}{ (1-\sig)^{(N+2)q-\frac{N+1}{2}}} \Big[\frac1{(1-\sig)^{N+1}}\biint_{Q(\rho, \theta)} u_+^\varep\,\dx\dt\Big]^{\frac12} \Big(\essup_{Q(R,\theta)} u\Big)^{\frac{2-\varep}{2}}\\
&\qquad\le \dl \essup_{Q(R,\theta)} u + \frac{\boldsymbol\gm^{\frac2\varep} \dl^{-\frac{2-\varep}{\varep}}\mathcal{A}_{\dl}^{\frac2\varep}}{(1-\sig)^{(N+2)q-\frac{N+1}{2}}} \Big[\frac1{(1-\sig)^{N+1}}\biint_{Q(\rho, \theta)} u_+^\varep\,\dx\dt\Big]^{\frac1\varep}.
\end{align*}
Plugging this estimate back to \eqref{Eq:A:4} and recalling the definition of $ \mathbf{I}$, we get
\begin{align}\nonumber
\essup_{Q(\sig\rho,\sig \theta)} u &\le \dl\Big[1+\Big(\frac{R}{\rho}\Big)^N\Big] \essup_{Q(R,\theta)} u + \frac12\Big(\frac{R}{R-\rho}\Big)^{N+2s} {\rm Tail} [u_+; Q(R, \theta)]\\ \label{Eq:A:5}
&\quad+  \frac{\boldsymbol\gm^{\frac2\varep} \dl^{-\frac{2-\varep}{\varep}}\mathcal{A}_{\dl}^{\frac2\varep}}{(1-\sig)^{(N+2)q-\frac{N+1}{2}}} \Big[\frac1{(1-\sig)^{N+1}}\biint_{Q(\rho, \theta)} u_+^\varep\,\dx\dt\Big]^{\frac1\varep},
\end{align}
 for all $Q(R,\theta)\subset E_T$, $\dl,\,\sig\in(0,1)$ and $\rho\in(0,R)$.

Next, we use \eqref{Eq:A:5} to perform an interpolation argument. To this end, introduce for  $n\in\nn_0$
\begin{equation*}
\left\{
\begin{array}{cc}
\dsty\rho_n=\rho - 2^{-n}(1-\sig)\rho,\quad\theta_n=\theta - 2^{-n}(1-\sig) \theta,\\ [5pt]
\dsty\widetilde{\rho}_n=\tfrac12(\rho_n +\rho_{n+1}),\quad \widetilde{\theta}_n=\tfrac12(\theta_n +\theta_{n+1}),\quad \sig_{n}=\frac{\rho_{n}}{\widetilde{\rho}_{n}}=\frac{\theta_{n}}{\widetilde{\theta}_{n}}.
\end{array}\right.
\end{equation*}
Here, we slightly abused the letter $\sig\in(0,1)$ as $\sig_n$ will play the role of $\sig$ in \eqref{Eq:A:5}.
Indeed, we apply \eqref{Eq:A:5} with $(\sig, \rho, R,  \theta)$ replaced by 
$(\sig_{n}, \widetilde{\rho}_n, \rho_{n+1}, \theta_{n+1})$. Notice also that $\theta_n\le\sig_n\theta_{n+1}$.
In this set-up, we obtain from \eqref{Eq:A:5} that for $n\in\nn_0$
\begin{align}\nonumber
\essup_{Q(\rho_n,\theta_n)} u &\le \dl\Big[1+\Big(\frac{\rho_{n+1}}{\widetilde{\rho}_n}\Big)^N\Big] \essup_{Q(\rho_{n+1},\theta_{n+1})} u\\ \nonumber
&\quad +  \frac12\Big(\frac{\rho_{n+1}}{\rho_{n+1}-\widetilde{\rho}_n}\Big)^{N+2s} {\rm Tail} [u_+; Q(\rho_{n+1}, \theta_{n+1})]\\ \nonumber
&\quad+\frac{\boldsymbol\gm^{\frac2\varep} \dl^{-\frac{2-\varep}{\varep}}\mathcal{A}_{\dl}^{\frac2\varep}}{(1-\sig_n)^{(N+2)q-\frac{N+1}{2}}}  \Big[\frac{1}{(1-\sig_n)^{N+1} \widetilde{\rho}_n^N \theta_{n+1} }\iint_{Q(\widetilde{\rho}_n, \theta_{n+1})} u_+^\varep\,\dx\dt\Big]^{\frac1\varep}\\\label{Eq:A:6}
&=:\mathbf{I}_1+\mathbf{I}_2+\mathbf{I}_3.
\end{align}
Note that $\mathcal{A}_{\dl}$ in \eqref{Eq:A:6} also depends on $\sig_n$. However, since $\sig_n\ge\sig$ and $\mathcal{A}_{\dl}(\sig_n)\le\mathcal{A}_{\dl}(\sig)$, such dependence has been replaced by $\sig$.
Let us continue to estimate the right-hand side of \eqref{Eq:A:6}. For the first term, we use the simple fact $\rho_{n+1}/\widetilde{\rho}_n\le2$ to bound it, such that
\[
\mathbf{I}_1\le \dl 4^{N}\essup_{Q(\rho_{n+1},\theta_{n+1})} u.
\] 
For the second term, we estimate
\begin{align*}
\mathbf{I}_2&=  \frac12\Big(\frac{2\rho_{n+1}}{\rho_{n+1}-\rho_n}\Big)^{N+2s} {\rm Tail} [u_+; Q(\rho_{n+1}, \theta_{n+1})] \\
 &\le \frac12 \Big(\frac{2^{n+2}}{1-\sig}\Big)^{N+2s} {\rm Tail} [u_+; Q(\sig \rho, \theta)].
\end{align*}
To estimate the third term, we use $\theta_{n+1}\ge\widetilde{\theta}_n$ to estimate
\begin{align*}
(1-\sig_n)^{N+1} \widetilde{\rho}_n^N \theta_{n+1}& 
=(\widetilde{\rho}_n-\rho_{n})^N(1-\sig_n)\theta_{n+1}\\
&= (\rho_{n+1}-\rho_n)^N(\widetilde{\theta}_{n}-\theta_n)\frac{\theta_{n+1}}{\widetilde{\theta}_n}\\
&\ge\tfrac12(\rho_{n+1}-\rho_n)^N(\theta_{n+1}-\theta_n) \\
&=\frac{\rho^N\theta (1-\sig)^{N+1}}{2^{1+(n+1)(N+1)}}.
\end{align*}
Similarly,
\[
1-\sig_n=1-\frac{\rho_n}{\widetilde{\rho}_n}=\frac{\rho_{n+1}-\rho_n}{\widetilde{\rho}_n}\ge\frac{1-\sig}{2^{n+2}}.
\]
Consequently, we get
\begin{align*}
\mathbf{I}_3\le\frac{\boldsymbol\gm^{\frac2\varep} \dl^{-\frac{2-\varep}{\varep}}\mathcal{A}_{\dl}^{\frac2\varep}\cdot 2^{\frac2\varep(n+2)(N+2)q}}{(1-\sig)^{q(N+2)-\frac{N+1}{2}}}  \Big[\frac1{(1-\sig)^{N+1}}\biint_{Q(\rho, \theta)} u_+^\varep\,\dx\dt\Big]^{\frac1\varep}.
\end{align*}
Here, the power of $2$ was $\frac1\varep[1+(n+1)(N+1)]+q(n+2)(N+2)-\frac12(N+1)$ and had been slightly simplified.
Combining all these estimates in \eqref{Eq:A:6} and denoting $\boldsymbol{b}:=2^{\frac{2q(N+2)}{\varep}}$ and $\bar{\dl}:=\dl 4^N$, we arrive at the recursive inequality for $n\in\nn_0$:
\begin{align*}
\essup_{Q(\rho_n,\theta_n)} u \le \bar{\dl} \essup_{Q(\rho_{n+1},\theta_{n+1})} u +\mathcal{B}_{\dl} \boldsymbol{b}^n
\end{align*}
where
\begin{align*}
\mathcal{B}_\dl&:=
 \frac12\Big(\frac{4}{1-\sig}\Big)^{N+2s} {\rm Tail} [u_+; Q(\sig \rho, 2\theta)]\\
&\qquad+\frac{\boldsymbol\gm^{\frac2\varep} \dl^{-\frac{2-\varep}{\varep}}\mathcal{A}_{\dl}^{\frac2\varep}\cdot 2^{\frac{4}{\varep} q(N+2)}}{(1-\sig)^{(N+2)q-\frac{N+1}{2}}}  \Big[\frac1{(1-\sig)^{N+1}}\biint_{Q(\rho, \theta)} u_+^\varep\,\dx\dt\Big]^{\frac1\varep}.
\end{align*}
Iterating this recursive inequality to obtain
\[
\essup_{Q(\rho_o,\theta_o)} u \le \bar{\dl}^n \essup_{Q(\rho_{n},\theta_{n})} u +\mathcal{B}_{\dl} \sum_{i=0}^{n-1}(\bar{\dl}\boldsymbol{b})^i.
\]
To conclude the proof, we choose $\bar\dl=1/(2\boldsymbol{b})$ and let $n\to\infty$.

Let us remark that the proof of the dimensionless version Theorem~\ref{Prop:A:2} follows the same lines with $\rho=1$.
\subsection{Further discussion}\label{S:discussion}
Previously, a standard procedure to obtain the nonlocal {\it elliptic} Harnack estimate hinges on, amongst other things, an interpolated boundedness estimate for sub-solutions of the following type:
\begin{equation}\label{Eq:sup-interpo-ell}
\essup_{K_{\frac12\rho}} u\le \dl {\rm Tail}[u_+; K_{\frac12\rho}] +\boldsymbol\gm(\dl)\Big(\bint_{K_\rho} u_+^2\,\dx\Big)^{\frac12}
\end{equation}
for any $\dl\in(0,1)$. With this one can circumvent the improved, weak Harnack estimate \eqref{Eq:intro:3} for super-solutions, and establish the Harnack estimate for solutions all the same.

A parabolic analogue of \eqref{Eq:sup-interpo-ell} with the tail \eqref{Eq:tail} would be
\begin{equation}\label{Eq:sup-interpo-par}
\essup_{Q_{\frac12\rho}} u\le \dl {\rm Tail}[u_+; Q(\tfrac12\rho,\rho^{2s})] +\boldsymbol\gm(\dl)\Big(\biint_{Q_\rho} u_+^2\,\dx\dt\Big)^{\frac12}.
\end{equation}
The purpose of this section is to explain that \eqref{Eq:sup-interpo-par} does not hold in general by providing explicit examples; compare with~\cite[Theorem~1.3]{Stromqvist-2} and note that \cite[Lemma~2.7]{Stromqvist-2} is unproven.

The following construction is similar to \cite[Example~5.2]{Kass-23}.
\begin{proposition}
Let $f(t)$ be an absolutely continuous, non-decreasing function in $[-1,1]$, and let $\lm$ be a positive parameter.
Consider the function
\begin{equation}\label{Eq:example}
v(x,t):=\lm f(t) + f^{\prime}(t)\chi_{K_3\setminus K_2}(x)\quad\text{for}\>(x,t)\in \rn\times[-1,1].
\end{equation}
There exists $\lm_o>0$ depending only on $N$, such that for any $\lm\in(0,\lm_o]$,  the so-defined  $v$ is a  weak sub-solution to $\pl_t u+(-\Dl)^s u=0$ in $K_1\times(-1,1]$ according to Definition~\ref{Def:solution}.
\end{proposition}
\begin{proof}
The notion of weak solution is standard; see~\cite{Liao-mod-24} for instance.
It is not hard to see that, such $v$ belongs to the function spaces in the definition of solution in \cite{Liao-mod-24}.
Hence, it remains to check the integral inequality
\begin{equation*} 
\begin{aligned}
	{\bf I}:=\int_{K_1} & v\vp \,\dx\bigg|_{t_1}^{t_2}
	-\int_{t_1}^{t_2}\int_{K_1}  v\pl_t\vp\, \dx\dt
	+\int_{t_1}^{t_2} \mathscr{E}\big(v(\cdot, t), \vp(\cdot, t)\big)\,\dt
	\le0,
\end{aligned}
\end{equation*}
where $[t_1,t_2]\subset [-1,1]$ and
\[
	 \mathscr{E}:=\int_{\rn}\int_{\rn}   \frac{ \big(v(x,t) - v(y,t)\big)\big(\vp(x,t) - \vp(y,t)\big)}{2|x-y|^{N+2s}}\,\dy\dx
\]
for all non-negative testing functions
\begin{equation*} 
\vp \in W^{1,2}\big([-1,1];L^2(K_1)\big)\cap L^2 \big([-1,1];W_o^{s,2}(K_1)\big).
\end{equation*}

First of all, since $f(t)$ and $\vp(\cdot, t)$ are absolutely continuous over $[-1,1]$, an integration by parts yields that
\begin{align*}
\int_{K_1}  v\vp \,\dx\bigg|_{t_1}^{t_2}
	-\int_{t_1}^{t_2}\int_{K_1}  v\pl_t\vp\, \dx\dt&=\int_{K_1}  \lm f \vp \,\dx\bigg|_{t_1}^{t_2}
	-\int_{t_1}^{t_2}\int_{K_1}  \lm f \pl_t\vp\, \dx\dt  \\
	&=\int_{t_1}^{t_2}\int_{K_1} \lm\vp f^{\prime}\, \dx\dt.
\end{align*}
Next, by symmetry, the nonlocal part gives for any fixed $t\in[-1,1]$ that
\begin{align*}
	 \mathscr{E}&=\int_{\rn}\int_{\rn}   \frac{ \big(v(x,t) - v(y,t)\big) \vp(x,t) }{|x-y|^{N+2s}}\,\dy\dx \\
	 &=\int_{\rn} \vp (-\Dl)^s v\,\dx=\int_{\rn} \vp f^{\prime} (-\Dl)^s \chi_{K_3\setminus K_2} \,\dx.
\end{align*}
The above calculation of $\mathscr{E}$ is legitimate because the last integral is convergent. In fact,
since ${\rm supp}\,\vp(\cdot, t)\subset K_1$, the effective integration domain of the last integral is actually $K_1$.
Furthermore, for $x\in K_1$ we can compute that
\begin{align*}
(-\Dl)^s \chi_{K_3\setminus K_2}(x)&=\int_{\rn}   \frac{  \chi_{K_3\setminus K_2}(x) - \chi_{K_3\setminus K_2}(y)   }{|x-y|^{N+2s}}\,\dy\\
&=-\int_{K_3\setminus K_2}   \frac{ \dy }{|x-y|^{N+2s}}\le -\frac{|K_3\setminus K_2|}{4^{N+2}}=:-\lm_o.
\end{align*}
Therefore, using these estimates and recalling $f^{\prime}\ge0$ due to its monotonicity, we have
\begin{align*}
{\bf I}&=\int_{t_1}^{t_2}\int_{K_1} \lm \vp f^{\prime}\, \dx\dt+ \int_{t_1}^{t_2}\int_{K_1} \vp f^{\prime} (-\Dl)^s \chi_{K_3\setminus K_2} \,\dx\dt\\
&\le\int_{t_1}^{t_2}\int_{K_1} (\lm-\lm_o)\vp f^{\prime}\, \dx\dt\le 0
\end{align*}
for any $\lm\in(0,\lm_o]$ as claimed.
\end{proof}

Suppose to the contrary that \eqref{Eq:sup-interpo-par} holds true for sub-solutions to the parabolic equation~\eqref{Eq:1:1} in $K_1\times(-1,1]$.
Consider such an estimate for the constructed sub-solution $v$ as in \eqref{Eq:example} with $f\equiv0$ in $[-1,0]$.
Let $\tau\in(0,1)$, and take the parabolic cylinder $Q=K_{\tau^{1/2s}}\times(0,\tau]$. Then \eqref{Eq:sup-interpo-par} implies that
\begin{equation}\label{Eq:sup-interpo-par-}
 v(0,\tau)\le \dl {\rm Tail}[v_+; K_{\frac12\tau^{1/2s}}\times(0,\tau]] +\boldsymbol\gm(\dl)\Big(\biint_{Q} v_+^2\,\dx\dt\Big)^{\frac12}. 
\end{equation}
Recalling the construction~\eqref{Eq:example}, the integral on the right-hand side of \eqref{Eq:sup-interpo-par-} becomes  
\[
\biint_{Q} v_+^2\,\dx\dt=\biint_{Q} [\lm f(t)]^2\,\dx\dt = \bint_{0}^{\tau} [\lm f(t)]^2\,\dt.
\]
Whereas the tail term in \eqref{Eq:sup-interpo-par-} is 
\begin{align*}
{\rm Tail}[v_+; K_{\frac12\tau^{1/2s}}\times(0,\tau]] &= \int_0^\tau\int_{ K_3 \setminus K_{2}}\frac{f^{\prime}(t)}{|x|^{N+2s}}\,\dx\dt+ \int_0^\tau\int_{\rn\setminus K_{\frac12\tau^{1/2s}}}\frac{\lm f(t)}{|x|^{N+2s}}\,\dx\dt \\
&= c_1 \int_0^\tau f^{\prime}(t)\,\dt+ c_2\bint_0^\tau \lm f(t)\,\dt \\
&= c_1 f(\tau)+ c_2\bint_0^\tau \lm f(t)\,\dt ,
\end{align*}
for some positive $c_1$ and $c_2$ depending on $N$ and $s$. The last integral can be further estimated by H\"older's inequality.
Plugging these estimates back to \eqref{Eq:sup-interpo-par-}, we obtain
\[
f(\tau)\le \dl c_1 f(\tau) + \big[\dl c_2+\boldsymbol\gm(\dl)\big] \Big(\bint_{0}^{\tau} [\lm f(t)]^2\,\dt\Big)^{\frac12}.
\]
Choosing $\dl \le \frac1{2c_1}$ we acquire that 
\[
f(\tau)\le   \Big[\tfrac{c_2}{c_1}+2\boldsymbol\gm(\tfrac1{2c_1})\Big] \Big(\bint_{0}^{\tau} [\lm f(t)]^2\,\dt\Big)^{\frac12}.
\]
To reach a contradiction, we select a sequence of functions $f_n(t)=t_+^n$ with $n\in\nn$. Plugging this selection into the last display, we obtain
\[
\sqrt{2n+1}\le \Big[\tfrac{c_2}{c_1}+2\boldsymbol\gm(\tfrac1{2c_1})\Big]\lm\quad\text{for any}\>n\in\nn.
\]
Thus, we have reached a contradiction.

\section{H\"older regularity}\label{S:Holder}

This section is devoted to the proof of Theorem~\ref{Thm:holder:1}. We follow the scheme devised in \cite{Liao-cvpd-24}; see also \cite{Byun-ampa,Liao-mod-24,NNSW} for further applications. The key element is the following expansion of positivity.
\begin{proposition}\label{Prop:expansion}
Let $u$ be a local, weak super-solution to \eqref{Eq:1:1}  in $E_T$.  Assume that $b$ satisfies \eqref{Eq:b:1} and that \eqref{Eq:global-thm} holds for some $\varep>0$.
 Suppose $u\ge0$ in $K_R(x_o)\times[t_o,t_o+R^{2s}]\subset E_T$ and for some constants $k>0$ and  $\al  \in(0,1]$, there holds
	\begin{equation*}
		\big|\big\{ u(\cdot, t_o) \ge k \big\}\cap K_{\varrho}(x_o) \big|
		\ge
		\al |K_\varrho |.
	\end{equation*}
There exist constants $\dl,\,\eta\in(0,1)$, $\widetilde{\boldsymbol\gm}>1$ depending only on the data $\{s,  N, C_o, C_1\}$ and $\al$,    such that 
either
\begin{equation*}
\Big(\frac{\rho}{R}\Big)^{\frac{2s\varep}{1+\varep}} \bigg[\bint_{t_o-R^{2s}}^{t_o}
 \Big(R^{2s}\int_{\rn\setminus K_R(x_o)}\frac{u_-(x,t)}{|x-x_o|^{N+2s}}\,\dx \Big)^{1+\varep}\dt\bigg]^{\frac1{1+\varep}}
>\eta k,
\end{equation*}
 or
\begin{equation*}
	u\ge\eta k
	\quad
	\mbox{a.e.~in $ K_{2\varrho}(x_o) \times\big( t_o+\tfrac12 \dl  \varrho^{2s},
	t_o+\dl  \varrho^{2s}\big],$}
\end{equation*}
provided  
\[
K_{4\rho}(x_o)\times\big(t_o, t_o+\dl  \varrho^{2s}\big]\subset K_R(x_o)\times[t_o,t_o+R^{2s}].
\]
Moreover,   $\dl\approx\al^{N+3}$ and $\eta\approx\al^q$ for some $q>1$ depending on the data $\{s, N, C_o, C_1\}$.
\end{proposition}

\noi Using the energy estimates in Proposition~\ref{Cor:1}, the proof of Proposition~\ref{Prop:expansion} follows in the similar way as \cite[Proposition~4.1]{Liao-cvpd-24}. The thrust is that in the range $s\in(\frac12,1)$, the nonlocal diffusion overtakes the drift-term with regard to the scaling of the equation and the energy estimates do not distinguish the drift or the drift-free case. Once the expansion of positivity is established, we can perform an induction argument as in \cite[\S~4]{Liao-cvpd-24} and conclude the proof of H\"older regularity.
We avoid repetition here and refer also to \cite{Byun-ampa, Liao-mod-24, NNSW} for a similar situation.
\section{Nonlocal elliptic weak Harnack estimates}\label{S:elliptic}
The main purpose of this section is to use our new method and show some new results in the elliptic case.
To this end, consider the fractional $p$-Laplacian type equation
\begin{equation}\label{Eq:p-lap}
{\rm p.v.}\int_{\rn} a(x,y)\frac{|u(x) - u(y)|^{p-2} (u(x) - u(y))}{|x-y|^{N+sp}}\,\dy=0\quad\text{weakly in}\>E,
\end{equation}
for some $p>1$, $s\in(0,1)$, and $a(x,y)$ is a symmetric, measurable function in $\rr^{2N}$ that satisfies $C_o\le a\le C_1$ for some positive $C_o,\,C_1$. The notion of solution is standard, cf.~\cite{Cozzi, Liao-cvpd-24}.
The tail is defined as
\begin{equation*}
{\rm Tail}[u; K_R(x_o)]:=  \Big(R^{sp}\int_{\rn\setminus K_R(x_o)}\frac{|u(x)|^{p-1}}{|x-x_o|^{N+sp}}\,\dx \Big)^{\frac1{p-1}}.
\end{equation*}

The following energy estimate is derived in \cite{Cozzi}; see also \cite[Corollary~3.2]{Liao-cvpd-24}. 
\begin{proposition}\label{Prop:p-energy}
Let $u$ be a local weak super-solution to \eqref{Eq:p-lap} in $E$.
There exists a constant $\boldsymbol \gm_* (C_o,C_1,p)>0$, such that
 	for all balls $K_r(x_o) \subset K_R(x_o)  \subset E$,
 	and every $k\in\rr$, 
\begin{align*}
	& \int_{K_r(x_o)}\int_{K_r(x_o)}  \frac{|w_-(x) - w_-(y)|^p}{|x-y|^{N+sp}}\,\dy\dx
	 +  \int_{K_r(x_o)}   \int_{\rn} \frac{w_{-}(x) w^{p-1}_+(y)}{|x-y|^{N+sp}}\,\dy\dx  \\
	 &\quad \le  \frac{ \boldsymbol \gm_* R^{(1-s)p}}{(R-r)^p} \int_{K_R(x_o)}  w^p_{-}\,\dx
	  +\frac{ \boldsymbol \gm_* R^{N}}{(R-r)^{N+sp}}\int_{K_R(x_o)} w_{-}\,\dx
	\big[{\rm Tail} [w_-; K_R(x_o) ]\big]^{p-1}, 
\end{align*}
where $w_-:=(u-k)_-$.
\end{proposition}

Once energy estimates are derived, we do not refer to equation~\eqref{Eq:p-lap} any more.
As in \cite{Cozzi}, we introduce the following {\it fractional DeGiorgi class}. 
\begin{definition}\label{Def:DG}
We call $u\in {\rm DG}^{s,p}_-(E; \boldsymbol \gm_*)$ if $u\in W_{\loc}^{s,p}(E)$, 
\[
\int_{\rn}\frac{|u(x)|^{p-1}}{1+|x|^{N+sp}}\,\dx<\infty,
\]
and satisfies the  estimate in Proposition~\ref{Prop:p-energy} for all $K_r(x_o) \subset K_R(x_o)  \subset E$ and every $k\in\rr$.
\end{definition}

\subsection{Main results}
Our first weak Harnack estimate examines the local integral contribution to the positivity of a function in ${\rm DG}^{s,p}_-(E; \boldsymbol \gm_*)$.
\begin{theorem}\label{Prop:WHI:1p}
 Let $u\in {\rm DG}^{s,p}_-(E; \boldsymbol \gm_*)$ be non-negative in $K_R(x_o)\subset E$. There exists a constant $\eta\in(0,1)$
 depending only on the data $\{s, p, N, \boldsymbol \gm_* \}$, such that
\begin{equation*}
\essinf_{K_{\frac12\rho}(x_o)} u + \Big(\frac{\rho}{R}\Big)^{\frac{sp}{p-1}}{\rm Tail}[u_-; K_R(x_o)] \ge 
\eta\Big(\bint_{K_\rho(x_o)} u^{p-1}\,\dx\Big)^{\frac1{p-1}},
\end{equation*}
provided that $K_{2\rho}(x_o)\subset K_R(x_o)$.
\end{theorem}

The second weak Harnack estimate measures the positivity contribution from the long-range behavior.
\begin{theorem}\label{Prop:WHI:2p}
 Let $u\in {\rm DG}^{s,p}_-(E; \boldsymbol \gm_*)$ be non-negative in $K_R(x_o)\subset E$. There exists a constant $\eta\in(0,1)$
 depending only on the data $\{s, p, N, \boldsymbol \gm_* \}$, such that
\begin{equation*}
\essinf_{K_{\frac12\rho}(x_o)} u + \Big(\frac{\rho}{R}\Big)^{\frac{sp}{p-1}}{\rm Tail}[u_-; K_R(x_o)] \ge 
 \eta{\rm Tail}[u_+; K_\rho(x_o)],
\end{equation*}
provided that $K_{2\rho}(x_o)\subset K_R(x_o)$.
\end{theorem}

\begin{remark}\upshape
One could adjust the estimates and replace $K_{\rho/2}(x_o)$ appearing under $\essinf$ by larger balls, for instance $K_{2\rho}(x_o)$, in both theorems.
\end{remark}
\subsection{What is new}\label{N-S-p}

In contrast to \cite[Proposition~6.8]{Cozzi}, the new feature of Theorem~\ref{Prop:WHI:1p} lies in that, the integral exponent has been improved to $p-1$. This improvement employs neither the Krylov-Safanov-type covering argument (see~\cite{DB-Trud}) nor a clustering lemma of DiBenedetto-Gianazza-Vespri-type (see~\cite[p.~16, Lemma~3.1]{DBGV-mono} and \cite{DIV-23}). However, like exhibited in the parabolic case, our approach essentially exploits the mixed truncation term in the energy estimate and captures information peculiar to the nonlocal structure. Needless to say, if one is allowed to use the equation \eqref{Eq:p-lap} again, the exponent in Theorem~\ref{Prop:WHI:1p} can be further raised. On the other hand, it is noteworthy that, in the weak Harnack estimate for {\it local} DeGiorgi classes, such an integral exponent has remained small and only been qualitatively determined, since the seminal work of DiBenedetto \& Trudinger, see~\cite{DB-Trud}. This directly impacted obtaining the conjectured sharp exponent $\frac1{p-1}$ in the Wiener-type criterion for quasi-minima; see~\cite{DBG-Wiener}. It is thus somewhat surprising that,  under the nonlocal framework, the integral exponent can be raised in a straightforward fashion.

Theorem~\ref{Prop:WHI:2p} is novel for the fractional DeGiorgi class. Surprisingly, it has been even overlooked for super-solutions to the fractional $p$-Laplace equation \eqref{Eq:p-lap}, cf.~\cite{DKP-2, Kass-09}.

A significance of DeGiorgi classes lies in that, it shows local properties of functions (e.g. H\"older regularity, Harnack's inequality, Wiener's criterion, etc.) hinges on minimizing functionals, rather than solving PDE's.
The fractional DeGiorgi class in Definition~\ref{Def:DG}  is not the most general one. Nevertheless, we choose this simplified version, in order to convey the main new ideas. With some technical complications, our approach can be generalized to various interesting settings.

\subsection{Basic elements}
The following DeGiorgi-type lemma is a time-independent version of \cite[Lemma~3.1]{Liao-cvpd-24}.
\begin{lemma}\label{Lm:DG:1p}
 Let $u\in {\rm DG}^{s,p}_-(E; \boldsymbol \gm_*)$ be non-negative in $K_R(x_o)\subset E$. 
 Let $k>0$  be a parameter.
There exists a  constant    $\nu\in(0,1)$ depending on the data $\{s, p, N, \boldsymbol \gm_* \}$, such that if
\begin{equation*}
	\big|\big\{
	u\le k\big\}\cap  K_{\varrho}(x_o)\big|
	\le
	\nu|K_{\varrho}|,
\end{equation*}
then either
\begin{equation*} 
\Big(\frac{\rho}{R}\Big)^{\frac{sp}{p-1}} 
{\rm Tail} [u_-; K_R(x_o) ]
>k,
\end{equation*}
or
\begin{equation*}
	u\ge \tfrac12k
	\quad
	\mbox{a.e.~in $ K_{\frac12\rho}(x_o)$,}
\end{equation*}
provided that $ K_{\rho}(x_o)\subset K_R(x_o)$.
\end{lemma}

The following measure shrinking lemma is a new element of this section. It   measures the smallness of the set $\{u\approx0\}$ according to a local integral.
\begin{lemma}\label{Lm:3:2p}
 Let $u\in {\rm DG}^{s,p}_-(E; \boldsymbol \gm_*)$ be non-negative in $K_R(x_o)\subset E$.
  Let $k>0$   be a parameter.
There exists a constant
 $\boldsymbol \gm>1$ depending only on the data $\{s, p, N, \boldsymbol \gm_*\}$, such that
  either 
  $$ 
\Big(\frac{\rho}{R}\Big)^{\frac{sp}{p-1}}  
{\rm Tail} [ u_-; K_R(x_o) ]
>k,
  $$ 
  or
\begin{equation*}
	\big|\big\{
	u\le  k \big\}\cap K_{\rho}(x_o)\big|
	\le  \frac{\boldsymbol\gm k^{p-1} }{ [u^{p-1}]_{K_{\rho}(x_o)}} |K_\rho|,
\end{equation*}
where $[\cdot]_{K_{\rho}(x_o)}$ is the integral average on $K_{\rho}(x_o)$,
provided that $K_{2\rho}(x_o)\subset K_R(x_o)$.
\end{lemma}
\begin{proof}
Assume $x_o$ is the origin.
Let us use the energy estimates of Proposition~\ref{Prop:p-energy} in $K_\rho\subset K_{2\rho}$. Namely,
\begin{align*}
& \int_{K_\rho}   \int_{\rn} \frac{w_{-}(x) w^{p-1}_+(y)}{|x-y|^{N+sp}}\,\dy\dx  \\
	 &\quad \le  \frac{ \boldsymbol \gm_* (2\rho)^{(1-s)p}}{\rho^p} \int_{K_{2\rho}}  w^p_{-}\,\dx
	  +\frac{ \boldsymbol \gm_* (2\rho)^{N}}{\rho^{N+sp}}\int_{K_{2\rho}} w_{-}\,\dx
	\big[{\rm Tail} [w_-; K_{2\rho} ]\big]^{p-1}\\
	&\quad\le \boldsymbol \gm\frac{k^{p}}{\rho^{sp}}  |K_\rho| +\boldsymbol \gm\frac{  k}{\rho^{sp}}|K_\rho|\big[{\rm Tail} [w_-; K_{2\rho} ]\big]^{p-1}.
\end{align*}
Here, we simply used $w_-=(u-k)_-\le k$ to obtain the last line.
The tail is estimated by
\begin{align*}
\big[{\rm Tail}[w_-; K_{2\rho}]\big]^{p-1} &= (2\rho)^{sp}\int_{\rn\setminus K_{2\rho}}\frac{w_-^{p-1}(x)}{|x|^{N+sp}}\,\dx\\
&= (2\rho)^{sp}\Big[\int_{K_R\setminus K_{2\rho}}\frac{w^{p-1}_-(x)}{|x|^{N+sp}}\,\dx+ \int_{\rn\setminus K_{R}}\frac{w^{p-1}_-(x)}{|x|^{N+sp}}\,\dx\Big]\\
&\le \boldsymbol \gm k^{p-1} + \boldsymbol \gm \Big(\frac{\rho}{R}\Big)^{sp}\big[{\rm Tail}[u_-; K_{R}]\big]^{p-1}\\
&\le \boldsymbol \gm k^{p-1}.
\end{align*}
In the last line, we enforced 
\[
\Big(\frac{\rho}{R}\Big)^{\frac{sp}{p-1}}  
{\rm Tail} [ u_-; K_R ]\le k.
\]
Hence, we have
\begin{equation}\label{Eq:shrink:0p}
 \int_{K_\rho}   \int_{\rn} \frac{w_{-}(x) w^{p-1}_+(y)}{|x-y|^{N+sp}}\,\dy\dx \le  \boldsymbol \gm\frac{k^{p}}{\rho^{sp}}  |K_\rho|.
\end{equation}

The left-hand side of the energy estimate \eqref{Eq:shrink:0p} is estimated by observing that $u_+^{p-1}\le c (u-k)_+^{p-1} + c k^{p-1}$ for some $c=c(p)$, and hence
\begin{align*}
 \int_{K_\rho}   \int_{\rn} \frac{w_{-}(x) w^{p-1}_+(y)}{|x-y|^{N+sp}}\,\dy\dx  &\ge \int_{K_\rho}   \int_{K_\rho} \frac{w_{-}(x) w^{p-1}_+(y)}{|x-y|^{N+sp}}\,\dy\dx\\
&\ge\int_{K_\rho}   \int_{K_\rho} \frac{w_{-}(x) [\frac1{c} u^{p-1}(y)-k^{p-1}]}{(2\rho)^{N+sp}}\,\dy\dx\\ 
&\ge\frac1{\boldsymbol \gm\rho^{sp}} \int_{K_\rho}w_{-}(x) \,\dx \bint_{K_\rho}  u^{p-1}(y)\,\dy - \boldsymbol \gm\frac{k^{p}}{\rho^{sp}}  |K_\rho| .
\end{align*}
Substituting these estimates back to the energy estimate, we arrive at 
\begin{align*}
\int_{K_\rho}w_{-}(x) \,\dx \cdot [u^{p-1}]_{K_\rho} \le \boldsymbol \gm k^{p}  |K_\rho| .
\end{align*}
The integral can be estimated from below by
\[
\int_{K_\rho}w_{-}(x) \,\dx \ge \tfrac12 k \big|\big\{u<\tfrac12 k\big\}\cap K_\rho\big|.
\]
Therefore, we have
\[
\big|\big\{u<\tfrac12 k\big\}\cap K_\rho\big|\le \frac{\boldsymbol \gm k^{p-1}}{[u^{p-1}]_{K_\rho}}|K_\rho|.
\]
Redefine $\tfrac12 k$ as $k$ and $2^{p-1} \boldsymbol\gm$ as $ \boldsymbol\gm$ to conclude.
\end{proof}

The next lemma shows the same spirit as the last one. However, the measure shrinking mechanism hinges upon a nonlocal integral.
\begin{lemma}\label{Lm:3:2p+}
 Let $u\in {\rm DG}^{s,p}_-(E; \boldsymbol \gm_*)$ be non-negative in $K_R(x_o)\subset E$.
 Let $k>0$ be a parameter. 
There exists a constant
 $\boldsymbol \gm>1$ depending only on the data $\{s, p, N, \boldsymbol \gm_* \}$,  such that
  either 
  $$  
  \Big(\frac{\rho}{R}\Big)^{\frac{sp}{p-1}}  
{\rm Tail} [ u_-;  K_R (x_o) ]
>k,
  $$ 
  or
\begin{equation*}
 \big|\big\{u \le  k \big\}\cap K_\rho(x_o)\big|\le  \frac{ \boldsymbol\gm k^{p-1}}{\big[{\rm Tail}[u_+; K_{\rho}(x_o)]\big]^{p-1}} |K_\rho|,
\end{equation*}
provided that $K_{2\rho}(x_o)\subset K_R(x_o)$.
\end{lemma}
\begin{proof}
Assume $x_o$ is the origin. The proof departs from \eqref{Eq:shrink:0p}.
The left-hand side of \eqref{Eq:shrink:0p} is estimated by observing that $u_+^{p-1}\le c (u-k)_+^{p-1} + c k^{p-1}$ for some $c=c(p)$, and 
 that when $|y|\ge \rho$ and $|x|\le \rho$, one has $|x-y|\le 2 |y|$.
As a result, we have
\begin{align*}
 \int_{K_\rho}   \int_{\rn} \frac{w_{-}(x) w^{p-1}_+(y)}{|x-y|^{N+sp}}\,\dy\dx  &\ge \int_{K_\rho}   \int_{\rn\setminus K_\rho} \frac{w_{-}(x) w^{p-1}_+(y)}{|x-y|^{N+sp}}\,\dy\dx\\
&\ge\int_{K_\rho}   \int_{\rn\setminus K_\rho} \frac{w_{-}(x) [\frac1{c} u_+^{p-1}(y)-k^{p-1}]}{(2|y|)^{N+sp}}\,\dy\dx\\ 
&\ge\frac1{\boldsymbol \gm\rho^{sp}} \int_{K_\rho}w_{-}(x) \,\dx \cdot \big[{\rm Tail}[u_+; K_{\rho}]\big]^{p-1} - \boldsymbol \gm\frac{k^{p}}{\rho^{sp}}  |K_\rho| .
\end{align*}
Substituting these estimates back to \eqref{Eq:shrink:0p}, we arrive at 
\begin{align*}
\int_{K_\rho}w_{-}(x) \,\dx \cdot \big[{\rm Tail}[u_+; K_{\rho}] \big]^{p-1} \le \boldsymbol \gm k^{p}  |K_\rho| .
\end{align*}
The left side integral is treated similarly as in the previous lemma. Therefore, we can conclude as before.
\end{proof}

\subsection{Proof of elliptic weak Harnack estimates}
Omit the reference to $x_o$ for simplicity.
Let $\nu$ be fixed in Lemma~\ref{Lm:DG:1p}, and choose $k_1$ and $k_2$ to satisfy
\[
 \frac{\boldsymbol\gm_1 k_1^{p-1} }{ [u^{p-1}]_{K_{\rho} }}\le\nu \qquad \text{and}\qquad  
\frac{ \boldsymbol\gm_2 k_2^{p-1}}{\big[{\rm Tail}[u_+; K_{\rho} ]\big]^{p-1}}\le \nu,
\]
where $\boldsymbol\gm_1$ and $\boldsymbol\gm_2$ are determined in Lemma~\ref{Lm:3:2p} and  Lemma~\ref{Lm:3:2p+} respectively.
According to Lemma~\ref{Lm:DG:1p} we have for $i=1,\,2$ that, either
\begin{equation*} 
\Big(\frac{\rho}{R}\Big)^{\frac{sp}{p-1}} 
{\rm Tail} [u_-; K_R ]
>k_i,
\end{equation*}
or
\begin{equation*}
	u\ge \tfrac12k_i
	\quad
	\mbox{a.e.~in $ K_{\frac12\rho}$.}
\end{equation*}
Therefore, the proof of Theorems~\ref{Prop:WHI:1p} \& \ref{Prop:WHI:2p} is finished.

%
%
%

\end{document}